\documentclass[11pt,a4paper]{article}

\usepackage{tikz}
\usepackage{amsthm}
\usepackage{amsmath}
\usepackage{amssymb}
\usepackage{mathrsfs}
\usepackage{geometry}
\usepackage{graphicx}
\usepackage{amsfonts}
\usepackage{epstopdf}
\usepackage{enumerate}
\usepackage{enumitem}
\usepackage{subcaption}
\usepackage{hyperref}

\usepackage{soul}

\usepackage[normalem]{ulem}

\numberwithin{equation}{section}

\allowdisplaybreaks

\newcommand{\Rmnum}[1]{\uppercase\expandafter{\romannumeral#1}} 

\def\Xint#1{\mathchoice
{\XXint\displaystyle\textstyle{#1}}%
{\XXint\textstyle\scriptstyle{#1}}%
{\XXint\scriptstyle\scriptscriptstyle{#1}}%
{\XXint\scriptscriptstyle\scriptscriptstyle{#1}}%
\!\int}
\def\XXint#1#2#3{{\setbox0=\hbox{$#1{#2#3}{\int}$ }
\vcenter{\hbox{$#2#3$ }}\kern-.6\wd0}}

\def\dashint{\Xint-}

\theoremstyle{plain}
\newtheorem{theorem}{Theorem}[section]
\newtheorem{proposition}[theorem]{Proposition}
\newtheorem{lemma}[theorem]{Lemma}
\newtheorem{corollary}[theorem]{Corollary}
\newtheorem{definition}[theorem]{Definition}

\theoremstyle{definition}
\newtheorem{remark}[theorem]{Remark}

\renewcommand{\thefootnote}{}

\AtEndDocument{
\bigskip{\footnotesize%
  \textsc{School of Mathematics, Hangzhou Normal University, Hangzhou 310036, China}  \par  
  \textit{E-mail address}: \texttt{gaojin@hznu.edu.cn} \par
}
\bigskip{\footnotesize%
  \textsc{Institut f\"ur Angewandte Mathematik, Universit\"at Bonn, Endenicher Allee 60, 53115 Bonn, Germany} \par
   \textsc{Current Address: Department of Mathematics, Aarhus University, 8000 Aarhus C, Denmark} \par
  \textit{E-mail address}: \texttt{yangmengqh@gmail.com}
}}

\begin{document}

\title{H\"older regularity of harmonic functions on metric measure spaces}
\author{Jin Gao and Meng Yang}
\date{}

\maketitle

\abstract{We introduce a H\"older regularity condition for harmonic functions on metric measure spaces and prove that, under a \emph{slow} volume regular condition and an upper heat kernel estimate, the H\"older regularity condition, the weak Bakry-\'Emery non-negative curvature condition, H\"older continuity of the heat kernel (with or without exponential terms), and the near-diagonal lower bound for the heat kernel are equivalent. As applications, first, we establish the validity of the so-called generalized reverse H\"older inequality on the Sierpi\'nski \emph{carpet} cable system, resolving an open problem left by Devyver, Russ, Yang (Int. Math. Res. Not. IMRN (2023), no. 18, 15537–15583). Second, we prove that two-sided heat kernel estimates \emph{alone} imply gradient estimates for the heat kernel on strongly recurrent fractal-like cable systems, improving the main results of the aforementioned paper. Third, we obtain H\"older (Lipschitz) estimates for the heat kernel on \emph{strongly recurrent} metric measure spaces, extending the classical Li-Yau gradient estimate for the heat kernel on Riemannian manifolds.}

\footnote{\textsl{Date}: \today}
\footnote{\textsl{MSC2020}: 28A80, 35K08}
\footnote{\textsl{Keywords}: H\"older regularity, gradient estimate, harmonic function, heat kernel, Sierpi\'nski carpet.}
\footnote{Jin Gao was supported by National Natural Science Foundation of China (Grant. No. 12271282), and Zhejiang Provincial Natural Science Foundation of China (Grant. No. LQN25A010019). Jin Gao acknowledges Jiaxin Hu and Eryan Hu for reading the preliminary draft. Part of the work was carried out while Meng Yang was participating a Mini-Symposia in \emph{the 9th European Congress of Mathematics (9ECM)} in Seville (Spain) between 15 and 19 July 2024, Meng Yang is very grateful to Fabrice Baudoin and Li Chen for the invitation to the symposia. The authors are grateful to Fabrice Baudoin for discussions on curvature conditions, and warmly thank Naotaka Kajino and Ryosuke Shimizu for pointing out an alternative proof of Proposition \ref{prop_main}. The authors thank the anonymous referees for their helpful comments and suggestions, which have improved the presentation of this paper.}

\renewcommand{\thefootnote}{\arabic{footnote}}
\setcounter{footnote}{0}

\section{Introduction}

Let us recall the following classical De Giorgi-Nash-Moser theorem from PDE theory. Let $D$ be a bounded domain in $\mathbb{R}^d$. Consider the following divergence form elliptic operator
$$Lu=\sum_{i,j=1}^d\frac{\partial}{\partial x_i}\left(a_{ij}\frac{\partial u}{\partial x_j}\right),$$
where $a_{ij}$, $i,j=1,\ldots,d$, are measurable functions in $D$ satisfying that $a_{ij}=a_{ji}$ for $i,j=1,\ldots,d$ and there exist some positive constants $\lambda, \Lambda$ such that
$$\lambda|\xi|^2\le\sum_{i,j=1}^da_{ij}(x)\xi_i\xi_j\le\Lambda|\xi|^2\text{ for any }x\in D,\xi=(\xi_1,\ldots,\xi_d)\in \mathbb{R}^d.$$
Then any weak solution $u\in W^{1,2}_{\mathrm{loc}}(D)$ of $Lu=0$ in $D$, that is, any ($L$-)harmonic function $u\in W^{1,2}_{\mathrm{loc}}(D)$, is locally H\"older continuous in $D$; equivalently, $u\in C^\alpha(D)$ for some $\alpha\in(0,1)$, where $\alpha$ is a H\"older exponent. One only needs the \emph{measurability} of the coefficients $a_{ij}$ to obtain \emph{H\"older regularity} of the weak solution $u$. To prove this theorem, very powerful iteration techniques were developed, including the De Giorgi iteration \cite{deGiorgi57} and the Moser iteration \cite{Moser61}, and were later applied in many areas; see \cite{CKS87,HS01,Cou03,GT12,BGK12,GHH24b} and the references therein for applications to analysis and heat kernel estimates on metric measure spaces. However, these iteration techniques give the existence of some H\"older exponent $\alpha$, which depends only on $d,\lambda, \Lambda$, but generally do not provide an explicit value of $\alpha$.

The main purpose of this paper is to consider H\"older regularity of harmonic functions on metric measure spaces, where an \emph{explicit} H\"older exponent can be obtained. These metric measure spaces arise from certain \emph{strongly recurrent} fractals.

\begin{figure}[ht]
\centering
\begin{minipage}[t]{0.48\textwidth}
\centering
\includegraphics[width=\textwidth]{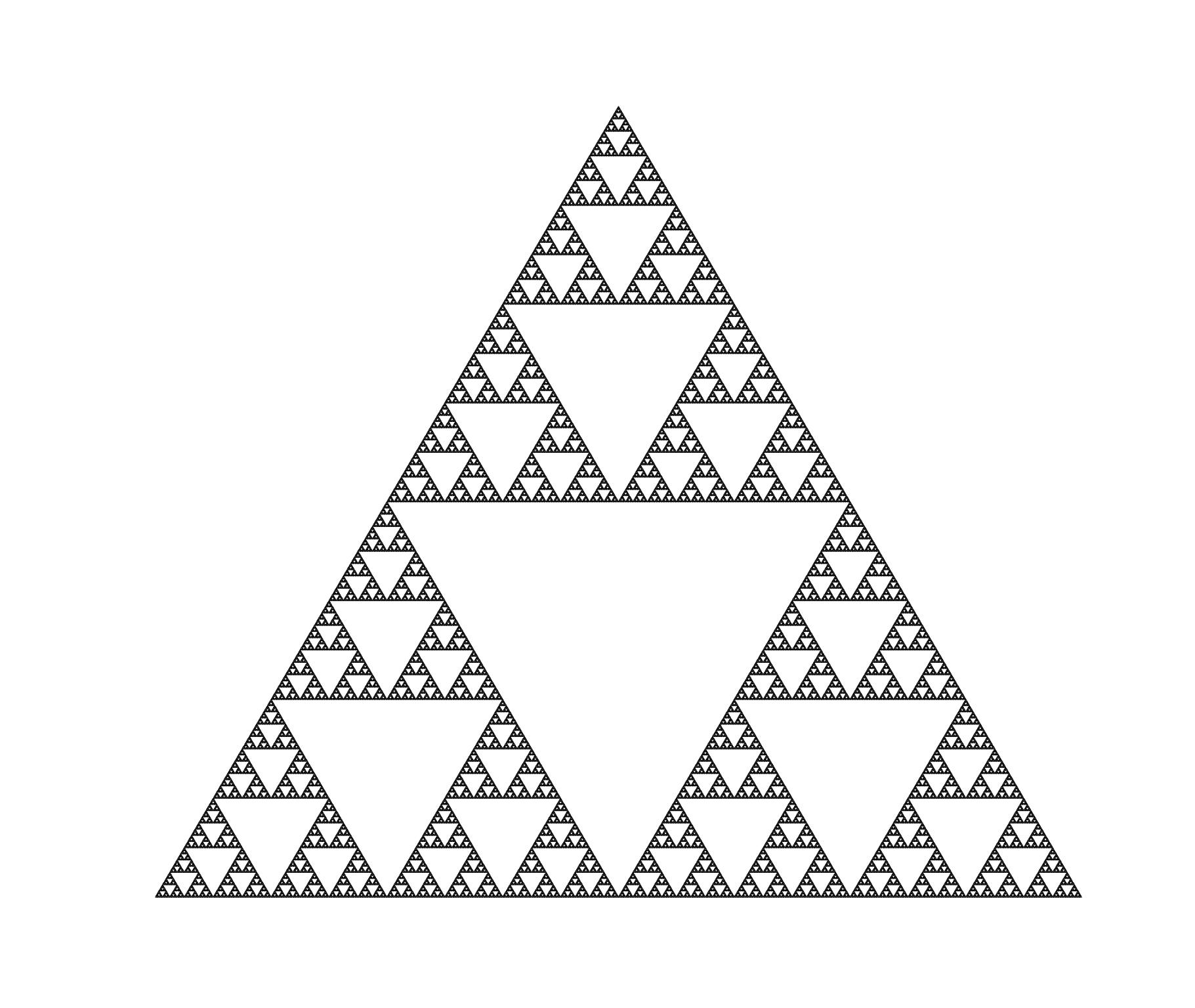}
\caption{The Sierpi\'nski gasket}\label{fig_SG}
\end{minipage}
\begin{minipage}[t]{0.48\textwidth}
\centering
\includegraphics[width=\textwidth]{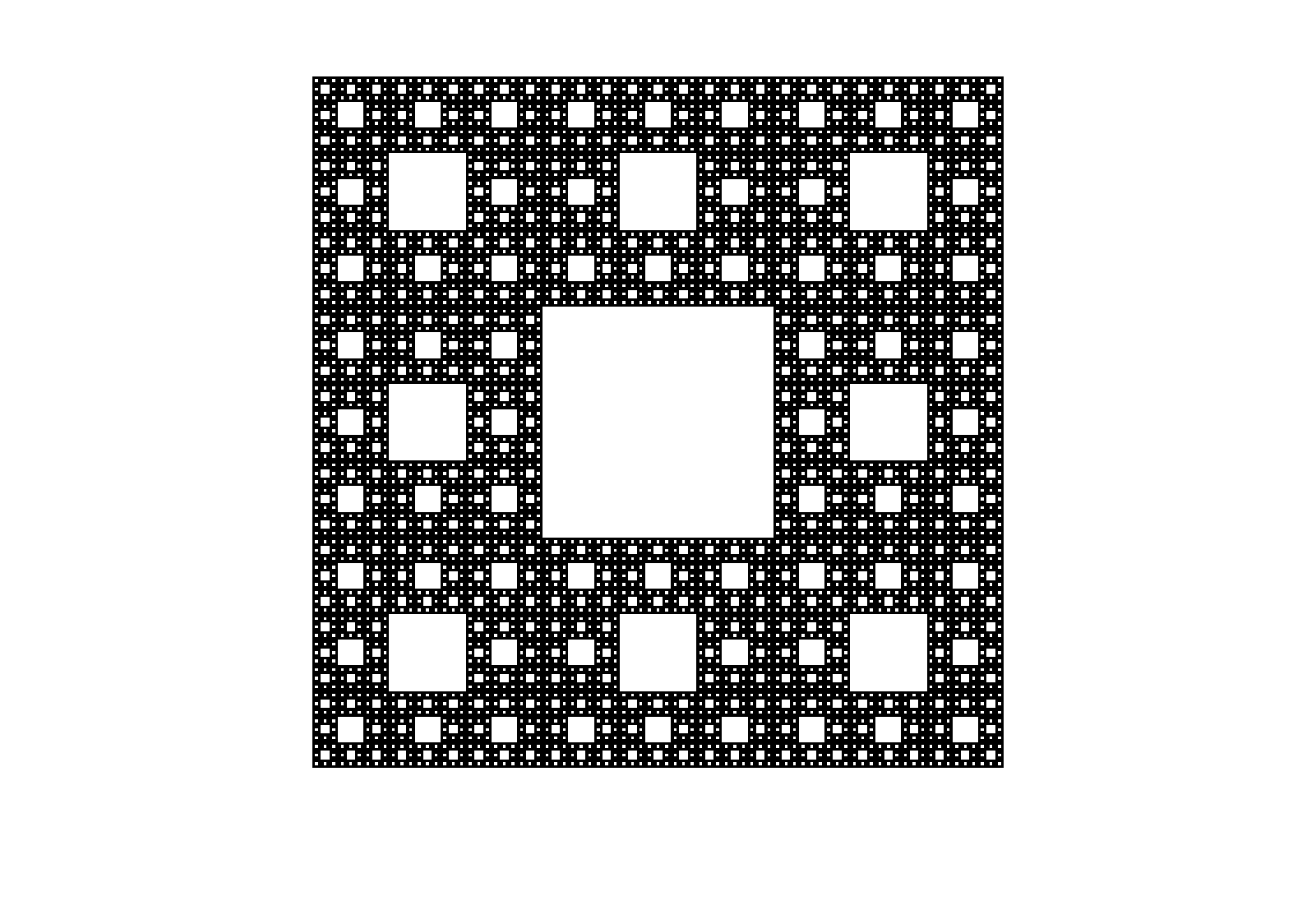}
\caption{The Sierpi\'nski carpet}\label{fig_SC}
\end{minipage}
\end{figure}

Fractals provide new examples exhibiting phenomena markedly different from those on Riemannian manifolds. On a complete non-compact Riemannian manifold, it was independently discovered by Grigor'yan \cite{Gri92} and Saloff-Coste \cite{Sal92,Sal95} that the following two-sided Gaussian estimate of the heat kernel
\begin{equation*}\label{eq_HK2}\tag*{HK(2)}
\frac{C_1}{V(x,\sqrt{t})}\exp\left(-C_2\frac{d(x,y)^2}{t}\right)\le p_t(x,y)\le\frac{C_3}{V(x,\sqrt{t})}\exp\left(-C_4\frac{d(x,y)^2}{t}\right),
\end{equation*}
is equivalent to the conjunction of the volume doubling condition and the scale-invariant $L^2$-Poincar\'e inequality on balls. However, one important estimate on fractals is the following two-sided \emph{sub-Gaussian} estimate of the heat kernel
\begin{equation*}\label{eq_HKbeta}\tag*{$\textrm{HK}(\beta)$}
\frac{C_1}{V(x,t^{1/\beta})}\exp\left(-C_2\left(\frac{d(x,y)}{t^{1/\beta}}\right)^{\frac{\beta}{\beta-1}}\right)\le p_t(x,y)\le\frac{C_3}{V(x,t^{1/\beta})}\exp\left(-C_4\left(\frac{d(x,y)}{t^{1/\beta}}\right)^{\frac{\beta}{\beta-1}}\right),
\end{equation*}
where $\beta$ is a new parameter called the walk dimension, which is typically strictly greater than 2 on fractals. Notably, for $\beta=2$, \ref{eq_HKbeta} is indeed \ref{eq_HK2}. For example, on the Sierpi\'nski gasket (see Figure \ref{fig_SG}), we have $\beta=\log5/\log2$, see \cite{BP88,Kig89}, while on the Sierpi\'nski carpet (see Figure \ref{fig_SC}), we have $\beta\approx2.09697$, see \cite{BB89,BB90,BBS90,BB92,HKKZ00,KZ92}.

Another purpose of this paper is to improve some recent results of Devyver, Russ and the second named author \cite{DRY23}, where gradient estimates for the heat kernel were considered. Gradient estimates for the heat kernel play an important role in the $L^p$-boundedness of the Riesz transform for $p>2$; see \cite{ACDH04,CJKS20}. To obtain pointwise upper \emph{sub-Gaussian} estimates for the gradient of the heat kernel on fractal-like cable systems, the authors of \cite{DRY23} introduced the so-called generalized reverse H\"older inequality and verified this condition on the Vicsek and Sierpi\'nski gasket cable systems, where the proof can be extended to cable systems corresponding to a class of p.c.f. self-similar sets. Since the harmonic extension algorithm is intrinsically needed, their proof does not apply to the Sierpi\'nski \emph{carpet} cable system. In this paper, we introduce a H\"older regularity condition, which is stronger than the generalized reverse H\"older inequality and, somewhat surprisingly, can be ensured by pointwise two-sided heat kernel estimates. As consequences, on the one hand, we establish the validity of the generalized reverse H\"older inequality on the Sierpi\'nski \emph{carpet} cable system, resolving an open problem left in \cite[Page 15545, Line 2]{DRY23}. On the other hand, to obtain gradient estimates for the heat kernel, we do not need to assume gradient-type conditions such as the generalized reverse H\"older inequality; instead, pointwise two-sided heat kernel estimates together with a slow volume regular condition suffice. This leads to an improvement of the main results of \cite[Theorem 1.1, Corollary 1.3]{DRY23}.

The next purpose of this paper is to generalize a recent result of Baudoin, Chen \cite[Theorem 3.3]{BC24} about Lipschitz estimates for the heat kernel on the unbounded Vicsek set. In this work, we obtain H\"older (Lipschitz) estimates for the heat kernel on strongly recurrent unbounded metric measure spaces, which may exhibit different behaviors at small and large scales. We follow a concise analytic approach, thereby avoiding the probabilistic framework of fractional diffusions introduced by Barlow \cite{Bar98}.

We conclude the Introduction by briefly stating our main result; a more detailed statement will be given in Section \ref{sec_state}.

Let $(X,d,m)$ be a metric measure space and let $(\mathcal{E},\mathcal{F})$ be a strongly local regular Dirichlet form on $L^2(X;m)$. A simplified consequence of our results is that, if $X$ is $\alpha$-regular for some $\alpha>0$, that is, for any metric ball $B(x,r)$, we have
$$V(x,r)=m(B(x,r))\asymp r^\alpha,$$
and if $(\mathcal{E},\mathcal{F})$ admits a heat kernel $p_t(x,y)$ satisfying \ref{eq_HKbeta} with walk dimension $\beta>\alpha$, then all harmonic functions associated with $(\mathcal{E},\mathcal{F})$ are $(\beta-\alpha)$-H\"older continuous. Moreover, the same H\"older regularity holds for the heat kernel. This result is motivated by classical examples on fractals, including the Sierpi\'nski gasket and the Sierpi\'nski carpet; see \cite{BP88,BB92}.

In fact, we consider more general and natural assumptions. We assume the volume regular condition
$$V(x,r)\asymp\Phi(r),$$
where
\begin{equation}\label{eq_Phi}
\Phi(r)=r^{\alpha_1}1_{(0,1)}(r)+r^{\alpha_2}1_{[1,+\infty)}(r),
\end{equation}
and $\alpha_1,\alpha_2>0$. \ref{eq_HKbeta} is also generalized by introducing a time-distance scaling function $\Psi$ of the form
\begin{equation}\label{eq_Psi}
\Psi(r)=r^{\beta_1}1_{(0,1)}(r)+r^{\beta_2}1_{[1,+\infty)}(r),
\end{equation}
which involves two walk dimensions: $\beta_1$ at small scales and $\beta_2$ at large scales. We further assume the strongly recurrent condition
\begin{equation}\label{eq_recurrent}
\alpha_i<\beta_i\text{ for }i=1,2.
\end{equation}
Under these assumptions, the H\"older regularity result is stated with two H\"older exponents: $\beta_1-\alpha_1$ at small scales and $\beta_2-\alpha_2$ at large scales.

More precisely, in an unbounded metric measure Dirichlet space, we introduce the H\"older regularity condition \ref{eq_HR} and prove that, assuming (\ref{eq_recurrent}), the volume regular condition \ref{eq_VPhi}, and the upper heat kernel estimate \ref{eq_UHK}, the following conditions are equivalent:
\begin{itemize}
\item the H\"older regularity condition \ref{eq_HR};
\item the weak Bakry-\'Emery non-negative curvature condition \ref{eq_wBE};
\item the H\"older continuity of the heat kernel \ref{eq_HHK};
\item the H\"older continuity of the heat kernel with exponential terms \ref{eq_HHKexp};
\item the near-diagonal lower bound for the heat kernel \ref{eq_NLE}.
\end{itemize}
Since we do not assume the existence of any gradient operator associated with the metric measure Dirichlet space, in particular, a ``\emph{carr\'e du champ}", the above equivalence generalizes the corresponding result of Coulhon, Jiang, Koskela, Sikora \cite[Theorem 1.2]{CJKS20}.

In addition, we extend a classical result from the setting of Riemannian manifolds to metric measure Dirichlet spaces. We prove that the conjunction of \ref{eq_UHK} and \ref{eq_HHKexp} implies \ref{eq_NLE} on metric measure Dirichlet spaces. This generalizes the fact that, on Riemannian manifolds, the conjunction of \emph{upper} heat kernel bounds and \emph{gradient} estimates implies a \emph{lower} heat kernel bound; see Cheeger, Yau \cite{CY81} and Li, Yau \cite{LY86}.

Throughout the paper, the letters $C,C_1,C_2,C_A,C_B$ will always refer to some positive constants and may change at each occurrence. The sign $\asymp$ means that the ratio of the two sides is bounded from above and below by positive constants. The sign $\lesssim$ ($\gtrsim$) means that the LHS is bounded by positive constant times the RHS from above (below).

\section{Statement of the main results}\label{sec_state}

Let $(X,d,m,\mathcal{E},\mathcal{F})$ be an unbounded metric measure Dirichlet (MMD) space, that is, $(X,d)$ is a locally compact separable unbounded metric space, $m$ is a positive Radon measure on $X$ with full support, and $(\mathcal{E},\mathcal{F})$ is a strongly local regular Dirichlet form on $L^2(X;m)$. Although the results of this paper also apply to bounded spaces, we focus on the unbounded case, where the space may exhibit different behaviors at small and large scales. Throughout this paper, we always assume that all metric balls are relatively compact, and that $(\mathcal{E},\mathcal{F})$ on $L^2(X;m)$ is conservative.

For any $x\in X$ and any $r\in(0,+\infty)$, denote the (metric) ball $B(x,r)=\{y\in X:d(x,y)<r\}$, denote $V(x,r)=m(B(x,r))$. If $B=B(x,r)$, then we denote $\delta B=B(x,\delta r)$ for any $\delta\in(0,+\infty)$. Let $C(X)$ denote the space of all real-valued continuous functions on $X$ and let $C_c(X)$ denote the space of all real-valued continuous functions on $X$ with compact support. Denote $\dashint_A=\frac{1}{m(A)}\int_A$ and $u_A=\dashint_Au\mathrm{d} m$ for any measurable set $A$ with $m(A)\in(0,+\infty)$ and any function $u$ such that the integral $\int_Au\mathrm{d} m$ is well-defined.

Consider the strongly local regular Dirichlet form $(\mathcal{E},\mathcal{F})$ on $L^2(X;m)$. Let $\Delta$ be the corresponding generator which is a non-negative definite self-adjoint operator. Let $\Gamma$ be the corresponding energy measure. Denote $\mathcal{E}_1(\cdot,\cdot)=\mathcal{E}(\cdot,\cdot)+(\cdot,\cdot)$, where $(\cdot,\cdot)$ is the inner product in $L^2(X;m)$. We refer to \cite{FOT11} for related results about Dirichlet forms.

We introduce several conditions in order to state our main results.

We say that the volume doubling condition \ref{eq_VD} holds if there exists $C_{VD}\in(0,+\infty)$ such that
\begin{equation*}\label{eq_VD}\tag*{VD}
V(x,2r)\le C_{VD}V(x,r)\text{ for any }x\in X,r\in(0,+\infty).
\end{equation*}

We say that the volume regular condition \ref{eq_VPhi} holds if there exists $C_{VR}\in(0,+\infty)$ such that
\begin{equation*}\label{eq_VPhi}\tag*{V($\Phi$)}
\frac{1}{C_{VR}}\Phi(r)\le V(x,r)\le C_{VR}\Phi(r)\text{ for any }x\in X,r\in(0,+\infty).
\end{equation*}

Consider the regular Dirichlet form $(\mathcal{E},\mathcal{F})$ on $L^2(X;m)$. Let $\left\{P_t\right\}$ be the corresponding heat semi-group. Let $\left\{X_t,t\ge0,\mathbb{P}_x,x\in X\backslash\mathcal{N}_0\right\}$ be the corresponding Hunt process, where $\mathcal{N}_0$ is a properly exceptional set, that is, $m(\mathcal{N}_0)=0$ and $\mathbb{P}_x(X_t\in\mathcal{N}_0\text{ for some }t>0)=0$ for any $x\in X\backslash\mathcal{N}_0$. For any bounded Borel function $f$, we have $P_tf(x)=\mathbb{E}_xf(X_t)$ for any $t>0$, for any $x\in X\backslash\mathcal{N}_0$.

The heat kernel $p_t(x,y)$ associated with the heat semi-group $\left\{P_t\right\}$ is a measurable function defined on $(0,+\infty)\times(X\backslash\mathcal{N}_0)\times(X\backslash\mathcal{N}_0)$ satisfying that:
\begin{itemize}
\item For any bounded Borel function $f$, for any $t>0$, for any $x\in X\backslash\mathcal{N}_0$, we have
$$P_tf(x)=\int_{X\backslash\mathcal{N}_0}p_t(x,y)f(y)m(\mathrm{d}y).$$
\item For any $t,s>0$, for any $x,y\in X\backslash\mathcal{N}_0$, we have
$$p_{t+s}(x,y)=\int_{X\backslash\mathcal{N}_0}p_t(x,z)p_s(z,y)m(\mathrm{d}z).$$
\item For any $t>0$, for any $x,y\in X\backslash\mathcal{N}_0$, we have $p_t(x,y)=p_t(y,x)$.
\end{itemize}
See \cite{GT12} for more details.

We say that the upper heat kernel estimate \ref{eq_UHK} holds if there exist a properly exceptional set $\mathcal{N}$ and $C_1$, $C_2\in(0,+\infty)$ such that for any $t\in(0,+\infty)$, for any $x,y\in X\backslash\mathcal{N}$, we have
\begin{equation*}\label{eq_UHK}\tag*{UHK($\Psi$)}
p_t(x,y)\le\frac{C_1}{V\left(x,\Psi^{-1}(t)\right)}\exp\left(-\Upsilon\left(C_2d(x,y),t\right)\right),
\end{equation*}
where
$$\Upsilon(R,t)=\sup_{s\in(0,+\infty)}\left(\frac{R}{s}-\frac{t}{\Psi(s)}\right)\asymp
\begin{cases}
\left(\frac{R}{t^{1/{\beta_1}}}\right)^{\frac{\beta_1}{\beta_1-1}},&\text{if }t<R,\\
\left(\frac{R}{t^{1/{\beta_2}}}\right)^{\frac{\beta_2}{\beta_2-1}},&\text{if }t\ge R.
\end{cases}$$
Then \ref{eq_UHK} can also be re-written as follows:
$$p_t(x,y)\le
\begin{cases}
\frac{C_1}{V\left(x,\Psi^{-1}(t)\right)}\exp\left(-C_2\left(\frac{d(x,y)}{t^{1/{\beta_1}}}\right)^{\frac{\beta_1}{\beta_1-1}}\right),&\text{if }t<d(x,y),\\
\frac{C_1}{V\left(x,\Psi^{-1}(t)\right)}\exp\left(-C_2\left(\frac{d(x,y)}{t^{1/{\beta_2}}}\right)^{\frac{\beta_2}{\beta_2-1}}\right),&\text{if }t\ge d(x,y).
\end{cases}
$$
If a corresponding lower bound, analogous to \ref{eq_UHK} but possibly with different constants $C_i$, also holds, then we say that \hypertarget{eq_HKPsi}{HK($\Psi$)} holds.

We say that the near-diagonal lower bound for the heat kernel \ref{eq_NLE} holds if there exist a properly exceptional set $\mathcal{N}$ and $C,\varepsilon\in(0,+\infty)$ such that for any $t\in(0,+\infty)$, for any $x,y\in X\backslash\mathcal{N}$ with $d(x,y)<\varepsilon\Psi^{-1}(t)$, we have
\begin{equation*}\label{eq_NLE}\tag*{NLE($\Psi$)}
p_t(x,y)\ge\frac{C}{V(x,\Psi^{-1}(t))}.
\end{equation*}
If the metric $d$ is furthermore assumed to satisfy the chain condition, then the conjunction of \ref{eq_UHK} and \ref{eq_NLE} is equivalent to \hyperlink{eq_HKPsi}{HK($\Psi$)}, see \cite[Section 6]{GT12}. Moreover, assuming \hyperlink{eq_HKPsi}{HK($\Psi$)}, $d$ satisfies the chain condition; see \cite[Theorem 2.11]{Mur20}. If, in addition, \ref{eq_VPhi} holds, then $2\le \beta_i\le \alpha_i+1$ for $i=1,2$, see \cite[THEOREM 2.1]{Mur24} and \cite[Corollary 4.9]{GHL03}.

Let $D$ be an open subset of $X$ and
$$\mathcal{F}_D=\text{ the }\mathcal{E}_1\text{-closure of }\mathcal{F}\cap C_c(D).$$
Let $f\in L^1_{\mathrm{loc}}(D)$. We say that $u\in\mathcal{F}$ is a solution of the Poisson equation $\Delta u=f$ in $D$ if
$$\mathcal{E}(u,\varphi)=\int_Df\varphi\mathrm{d}m\text{ for any }\varphi\in\mathcal{F}\cap C_c(D).$$
If $\Delta u=f$ in $D$ with $f\in L^2(D)$, then the above equation also holds for any $\varphi\in\mathcal{F}_D$. We say that $u\in\mathcal{F}$ is harmonic in $D$ if $\Delta u=0$ in $D$.

Denote $\gamma_i=\beta_i-\alpha_i\in(0,+\infty)$ for $i=1,2$, then
$$\left(\frac{\Psi}{\Phi}\right)(r)=r^{\gamma_1}1_{(0,1)}(r)+r^{\gamma_2}1_{[1,+\infty)}(r).$$

We say that the H\"older regularity condition \ref{eq_HR} holds if there exists $C_H\in(0,+\infty)$ such that for any ball $B=B(x_0,r)$, for any $u\in\mathcal{F}$ which is harmonic in $2B$, for $m$-a.e. $x,y\in B$ with $x\ne y$, we have
\begin{equation*}\label{eq_HR}\tag*{HR($\Phi,\Psi$)}
\frac{\Phi(d(x,y))}{\Psi(d(x,y))}|u(x)-u(y)|\le C_H\frac{\Phi(r)}{\Psi(r)}\dashint_{2B}|u|\mathrm{d}m.
\end{equation*}
The above condition can be viewed as a generalization of the so-called generalized reversed H\"older inequality \ref{eq_GRH} introduced in \cite{DRY23}.

We say that the H\"older estimate for the heat kernel \ref{eq_HHK} holds if there exists $C_{HHK}\in(0,+\infty)$ such that for any $t\in(0,+\infty)$, for any $x_1,x_2,y_1,y_2\in X$, we have
\begin{equation*}\label{eq_HHK}\tag*{HHK($\Phi,\Psi$)}
|p_t(x_1,y_1)-p_t(x_2,y_2)|\le C_{HHK}\frac{\left(\frac{\Psi}{\Phi}\right)(d(x_1,x_2))+\left(\frac{\Psi}{\Phi}\right)(d(y_1,y_2))}{t}.
\end{equation*}
The above condition implies that for any $t\in(0,+\infty)$, for any $x\in X$, $p_t(x,\cdot)$ is $\gamma_1$-H\"older continuous at small scales and $\gamma_2$-H\"older continuous at large scales.

We say that the H\"older estimate with exponential terms for the heat kernel \ref{eq_HHKexp} holds if there exist $C_1,C_2\in(0,+\infty)$ such that for any $t\in(0,+\infty)$, for any $x,y_1,y_2\in X$, we have
\begin{align*}
&|p_t(x,y_1)-p_t(x,y_2)|\nonumber\\
&\le C_1\frac{\left(\frac{\Psi}{\Phi}\right)(d(y_1,y_2))}{t}\left(\exp\left(-\Upsilon\left(C_2d(x,y_1),t\right)\right)+\exp\left(-\Upsilon\left(C_2d(x,y_2),t\right)\right)\right).\label{eq_HHKexp}\tag*{$\text{HHK}_{\text{exp}}(\Phi,\Psi)$}
\end{align*}

We say that the weak Bakry-\'Emery non-negative curvature condition \ref{eq_wBE} holds if there exists $C_{wBE}\in(0,+\infty)$ such that for any $t\in(0,+\infty)$, for any $f\in L^\infty(X;m)$, for any $x,y\in X$, we have
\begin{equation*}\label{eq_wBE}\tag*{wBE($\Phi,\Psi$)}
\frac{\Phi(d(x,y))}{\Psi(d(x,y))}|P_tf(x)-P_tf(y)|\le C_{wBE}\frac{\Phi(\Psi^{-1}(t))}{\Psi(\Psi^{-1}(t))}\lVert {f}\rVert_{L^\infty(X;m)}.
\end{equation*}
The above condition can be viewed as a generalization of the weak Bakry-\'Emery non-negative curvature condition introduced in \cite[Definition 3.1]{ABCRST3}. Roughly speaking, letting $t\to0+$ in \ref{eq_wBE} suggests that the space is strongly recurrent.

The main result of this paper is as follows.

\begin{theorem}\label{thm_main}
Let $(X,d,m,\mathcal{E},\mathcal{F})$ be an unbounded MMD space satisfying \ref{eq_VPhi} and\\
\noindent \ref{eq_UHK}, where $\Phi$ and $\Psi$ are given by (\ref{eq_Phi}) and (\ref{eq_Psi}), respectively, with parameters $\alpha_i,\beta_i$, $i=1,2$, satisfying the strongly recurrent condition (\ref{eq_recurrent}). Then the followings are equivalent.
\begin{enumerate}[label=(\arabic*)]
\item \ref{eq_HR}.
\item \ref{eq_wBE}.
\item \ref{eq_HHK}.
\item \ref{eq_HHKexp}.
\item \ref{eq_NLE}.
\end{enumerate}
\end{theorem}

The key novelty of the above result is that \emph{pointwise} two-sided heat kernel estimates suffice to give \emph{explicit} H\"older regularity for harmonic functions. In the existing literature (see, for example, \cite[Section 5.3]{HS01}, \cite[Proposition 4.5]{Cou03}, \cite[Section 5.3]{GT12}, and \cite[Section 4.2.1]{BGK12}), approaches based on the Moser or De Giorgi iteration, together with oscillation inequalities, typically give only H\"older continuity with an exponent that depends \emph{implicitly} on other constants.

Our proof will follow the implication diagrams below:
$$\text{\ref{eq_HR}}\Rightarrow\text{\ref{eq_wBE}}\Rightarrow\text{\ref{eq_HHK}}\Rightarrow\text{\ref{eq_NLE}}\Rightarrow\text{\ref{eq_HR}},$$
and
$$\text{\ref{eq_HR}}\Rightarrow\text{\ref{eq_HHKexp}}\Rightarrow\text{\ref{eq_HHK}}.$$
The implication ``\ref{eq_HHKexp}$\Rightarrow$\ref{eq_HHK}" is trivial, since it follows by simply neglecting the exponential terms. The proofs of ``\ref{eq_wBE}$\Rightarrow$\ref{eq_HHK}$\Rightarrow$\ref{eq_NLE}" are standard and will be given for completeness. The proofs of ``\ref{eq_HR}$\Rightarrow$\ref{eq_wBE}" and ``\ref{eq_HR}$\Rightarrow$\ref{eq_HHKexp}" will rely on techniques similar to the gradient estimates developed in \cite{DRY23,CJKS20}; see also \cite{Jiang15, JKY14}. The proof of ``\ref{eq_NLE}$\Rightarrow$\ref{eq_HR}" is new and will use a recent technique of resistance estimates in \cite[Theorem 6.27]{KS24a}. An alternative approach suggested by Kajino and Shimizu will also be provided, which relies on Kigami’s theory of resistance forms \cite{Kig12}.

\begin{remark}
In the setting of fractional metric spaces with fractional diffusions introduced in \cite[Section 3]{Bar98}, when the dimension $\alpha$ of the state space is strictly smaller than the walk dimension $\beta$, that is, when $\alpha_1=\alpha_2=\alpha<\beta_1=\beta_2=\beta$ in our notation, it was proved in \cite[Lemma 3.4, Remark 3.6, Theorem 3.7]{ABCRST3} that
$$\text{\ref{eq_wBE}}\Rightarrow\text{\ref{eq_HHK}}\Rightarrow\text{\ref{eq_NLE}}\Rightarrow\text{\ref{eq_wBE}}.$$
To prove ``\ref{eq_NLE}$\Rightarrow$\ref{eq_wBE}", the H\"older regularity result for resolvents \cite[Theorem 3.40]{Bar98} obtained within a probabilistic framework was intrinsically needed.
\end{remark}

\begin{remark}
In the fractal setting, \ref{eq_HHK} was also established on the unbounded Sierpi\'nski gasket (see \cite[Theorem 1.5 (c)]{BP88}) and on the unbounded Sierpi\'nski carpet (see \cite[Theorem 1.1 (e)]{BB92}).
\end{remark}

\begin{remark}\label{rmk_GLY}
The above equivalence can be regarded as a H\"older-type analogue of the gradient-type result in \cite[Theorem 1.2]{CJKS20}. In that work, it was proved that on a non-compact doubling Dirichlet metric measure space $(X,d,\mu,\mathscr{E})$ endowed with a ``carr\'e du champ", assuming the upper Gaussian bound $(UE)$ for the heat kernel and a local $L^\infty$-Poincar\'e inequality $(P_{\infty,\text{loc}})$, the following conditions are equivalent.
\begin{itemize}
\item \ref{eq_RH} Quantitative reverse $L^\infty$-H\"older inequality for gradients of harmonic functions: There exists $C\in(0,+\infty)$ such that for any ball $B=B(x_0,r)$, for any $u$ which is harmonic in $2B$, we have
\begin{equation*}\label{eq_RH}\tag*{($\text{RH}_\infty$)}
\lVert {|\nabla u|}\rVert_{L^\infty(B)}\le \frac{C}{r}\dashint_{2B}|u|\mathrm{d}\mu.
\end{equation*}
\item \ref{eq_GLY} Pointwise Li-Yau gradient estimate for the heat kernel: There exist $C,c\in(0,+\infty)$ such that
\begin{equation*}\label{eq_GLY}\tag*{($\text{GLY}_\infty$)}
|\nabla_xh_t(x,y)|\le\frac{C}{\sqrt{t}V(y,\sqrt{t})}\exp\left\{-c\frac{d(x,y)^2}{t}\right\}
\end{equation*}
for any $t\in(0,+\infty)$, for $\mu$-a.e. $x,y\in X$.
\item \ref{eq_Ginfty} $L^\infty$-boundedness of the gradient of the heat semi-group $|\nabla H_t|$: There exists $C\in(0,+\infty)$ such that for any $t\in(0,+\infty)$, we have
\begin{equation*}
\lVert {|\nabla H_t|}\rVert_{\infty\to\infty}\le\frac{C}{\sqrt{t}},
\end{equation*}
that is, for any $f\in L^\infty(X)$, we have
\begin{equation*}\label{eq_Ginfty}\tag*{($\text{G}_\infty$)}
\lVert {|\nabla H_tf|}\rVert_{L^\infty(X)}\le\frac{C}{\sqrt{t}}\lVert {f}\rVert_{L^\infty(X)}.
\end{equation*}
\item \ref{eq_GBE} A generalized Bakry-\'Emery condition: There exist $C,c\in(0,+\infty)$ such that for any $f\in W^{1,2}(X)$, for any $t\in(0,+\infty)$ and $\mu$-a.e. $x\in X$, we have
\begin{equation*}\label{eq_GBE}\tag*{($\text{GBE}$)}
\lvert \nabla H_tf(x)\rvert^2\le C H_{ct} (\lvert \nabla f\rvert^2)(x).
\end{equation*}
\end{itemize}

Our conditions can be viewed as H\"older-type analogues of the above gradient conditions, as illustrated in Figure \ref{fig_holdergrad}.

\begin{figure}[ht]
\centering
\begin{tabular}{|c|c|}
\hline
H\"older-type conditions& gradient-type conditions\\
\hline
\ref{eq_HR}&\ref{eq_RH}\\
&\\
\ref{eq_HHKexp}&\ref{eq_GLY}\\
&\\
\ref{eq_wBE}&\ref{eq_Ginfty}\\
\hline
\end{tabular}
\caption{The correspondence between H\"older-type and gradient-type conditions}\label{fig_holdergrad}
\end{figure}

In their proof of \ref{eq_RH}, they established a reproducing formula for harmonic functions using the finite propagation speed property, which may not be available on general MMD spaces. We leave it as an open question to extend the reproducing formula technique to prove \ref{eq_HR} on MMD spaces or the following \ref{eq_GRH} on cable systems, and in particular to develop a corresponding finite propagation speed theory in these settings.
\end{remark}

\begin{remark}
Let us compare the stability of the gradient-type and H\"older-type conditions as follows.

By a standard argument, under $(UE)$, any of the equivalent gradient-type conditions implies the near-diagonal Gaussian lower bound $(NLE)$ for the heat kernel. In particular, assuming $(UE)$, we have \ref{eq_GLY} implies $(NLE)$. However, the converse is usually not true. Indeed, the conjunction of $(UE)$ and $(NLE)$, or equivalently, \ref{eq_HK2}, is equivalent to the conjunction of the volume doubling condition and the scale-invariant $L^2$-Poincar\'e inequality on balls, and is therefore quasi-isometry \emph{stable}. By contrast, to obtain \ref{eq_GLY}, some curvature assumptions are usually needed, see \cite{LY86,BG13}, and hence \ref{eq_GLY} is quasi-isometry \emph{unstable}.

For comparison, by our main result (Theorem \ref{thm_main}), under \ref{eq_UHK}, all the H\"older-type conditions are equivalent to \ref{eq_NLE}. By the quasi-isometry stability of \ref{eq_UHK} (see \cite{AB15}), as well as the quasi-isometry stability of the conjunction of \ref{eq_UHK} and \ref{eq_NLE} (see \cite{GHL15}), our H\"older-type conditions are quasi-isometry \emph{stable}. Roughly speaking, this stability stems from the strongly recurrent condition, which forces our MMD spaces to behave analogously to $\mathbb{R}^1$ in the classical Riemannian setting.
\end{remark}

Let $(V,E)$ be an infinite, locally bounded, connected (undirected) graph. Let\\
\noindent $(X,d,m,\mathcal{E},\mathcal{F})$ be the corresponding unbounded cable system, equipped with the elementary gradient operator $\nabla$, see \cite[Section 3]{DRY23} for more details. Let $\alpha_1=1$, $\beta_1=2$, $\alpha_2=\alpha$, $\beta_2=\beta$.

We say that the generalized reverse H\"older inequality \ref{eq_GRH} holds if there exists $C_H\in(0,+\infty)$ such that for any ball $B=B(x_0,r)$, for any $u\in\mathcal{F}$ which is harmonic in $2B$, we have
\begin{equation}\label{eq_GRH}\tag*{GRH($\Phi,\Psi$)}
\lVert {\lvert\nabla u\rvert}\rVert_{L^\infty(B)}\le C_H\frac{\Phi(r)}{\Psi(r)}\dashint_{2B}|u|\mathrm{d}m.
\end{equation}
It is obvious that, on an unbounded cable system, \ref{eq_HR} implies \ref{eq_GRH}. Indeed, assuming \ref{eq_HR}, the LHS of \ref{eq_HR} reduces to $|u(x)-u(y)|/d(x,y)$ for $d(x,y)\in(0,1)$. Letting $y\to x$ and using the Lebesgue density theorem, we obtain the limit $|\nabla u(x)|$. Taking the essential supremum over $x\in B$, we have \ref{eq_GRH}. Roughly speaking, \ref{eq_HR} encodes more H\"older continuity information than \ref{eq_GRH} about the behavior of harmonic functions at \emph{large} scales.

We give some results as important applications of our main result, Theorem \ref{thm_main}, to unbounded cable systems. First, since two-sided heat kernel estimates have already been established for a large class of fractals (see \cite{Bar98,HK99,Kig01book,Kig09}), including the Sierpi\'nski carpet (see \cite{BB92,BB99,HKKZ00}), we have \ref{eq_GRH} holds on a large class of fractal-like cable systems, including the Sierpi\'nski \emph{carpet} cable system, a case whose validity was left open in \cite{DRY23}.

\begin{corollary}
Let $(X,d,m,\mathcal{E},\mathcal{F})$ be an unbounded cable system satisfying \ref{eq_VPhi} and\\
\noindent \hyperlink{eq_HKPsi}{HK($\Psi$)} with $\alpha_1=1$, $\beta_1=2$, and $\alpha_2=\alpha<\beta=\beta_2$. Then \ref{eq_GRH} holds. In particular, \ref{eq_GRH} holds on the Sierpi\'nski \emph{carpet} cable system with $\alpha=\log8/\log3<2<\beta\approx2.09697$.
\end{corollary}

\begin{remark}
In \cite{DRY23}, the authors verified \ref{eq_GRH} \emph{directly} on the Vicsek cable system and the Sierpi\'nski gasket cable system, using \emph{certain} H\"older regularity\footnote{The H\"older regularity follows indeed from a certain extension algorithm for harmonic functions.} of harmonic functions on the corresponding fractals \cite{Str00}, which is also available on some class of p.c.f. self-similar sets, but not available on general fractals, especially on the Sierpi\'nski carpet.
\end{remark}

Second, we improve the results about gradient estimates for the heat kernel as follows.

\begin{corollary}\label{cor_cable_grad}
Let $(X,d,m,\mathcal{E},\mathcal{F})$ be an unbounded cable system satisfying \ref{eq_VPhi} and\\
\noindent \hyperlink{eq_HKPsi}{HK($\Psi$)} with $\alpha_1=1$, $\beta_1=2$, and $\alpha_2=\alpha<\beta=\beta_2$.
\begin{enumerate}[label=(\arabic*)]
\item Then \ref{eq_HHKexp} holds. In particular, we have the gradient estimate \ref{eq_GHK} for the heat kernel as follows: there exist $C_1,C_2\in(0,+\infty)$ such that for any $t\in(0,+\infty)$, for $m$-a.e. $x,y\in X$, we have
\begin{align*}\label{eq_GHK}\tag*{GHK($\Phi,\Psi$)}
|\nabla_yp_t(x,y)|&\le\frac{C_1\Phi(\Psi^{-1}(t))}{tV\left(x,\Psi^{-1}(t)\right)}\exp\left(-\Upsilon\left(C_2d(x,y),t\right)\right).
\end{align*}
\item There exist $C_3,C_4\in(0,+\infty)$ such that for any $t\in(0,+\infty)$, for $m$-a.e. $x,y\in X$, we have
$$\left\vert \nabla_yp_t(x,y)\right\vert\le C_3\frac {\Phi(\Psi^{-1}(t))}{t}p_{C_4t}(x,y).$$
\end{enumerate}
\end{corollary}

\begin{remark}
(1) constitutes an improvement of \cite[Theorem 1.1]{DRY23}, while (2) improves upon \cite[Corollary 1.3]{DRY23}, where \ref{eq_GRH} was assumed. In the present work, we obtain the gradient estimate \ref{eq_GHK} \emph{without} imposing gradient-type conditions such as \ref{eq_GRH}.
\end{remark}

As a direct consequence of \ref{eq_HHKexp}, we obtain the following H\"older estimate for the heat kernel on unbounded MMD spaces.

\begin{corollary}\label{cor_MMD_holder}
Let $(X,d,m,\mathcal{E},\mathcal{F})$ be an unbounded MMD space satisfying \ref{eq_VPhi}\\
\noindent and \hyperlink{eq_HKPsi}{HK($\Psi$)}. Then there exist $C,c\in(0,+\infty)$ such that for any $t\in(0,+\infty)$, for any $x,y_1,y_2\in X$, we have
$$|p_t(x,y_1)-p_t(x,y_2)|\le C\frac{\left(\frac{\Psi}{\Phi}\right)(d(y_1,y_2))}{\left(\frac{\Psi}{\Phi}\right)(\Psi^{-1}(t))}\left(p_{ct}(x,y_1)+p_{ct}(x,y_2)\right).$$
\end{corollary}

\begin{remark}
If $\beta_1=\alpha_1+1$, then the above H\"older estimate gives a local Lipschitz estimate, which will imply the existence and some estimate of the gradient of the heat kernel. We give some examples as follows.
\begin{itemize}
\item Unbounded cable systems satisfying \ref{eq_VPhi} and \hyperlink{eq_HKPsi}{HK($\Psi$)}. In this case, $\alpha_1=1$ and $\beta_1=2$. This indeed recovers Corollary \ref{cor_cable_grad}.
\item The unbounded Vicsek set. Here $\alpha_1=\alpha_2=\log5/\log3$ and $\beta_1=\beta_2=\log15/\log3$, hence the above H\"older estimate is indeed a Lipschitz estimate, which was also given by Baudoin, Chen \cite[Theorem 3.3]{BC24}. However, a resolvent kernel estimate \cite[Theorem 3.40]{Bar98} was intrinsically needed in their proof.
\item Blowups of the Vicsek set by another fractal. In this case, $\alpha_1=\log5/\log3$ and $\beta_1=\log15/\log3$, while $\alpha_2$ and $\beta_2$ will be given by the Hausdorff dimension and the walk dimension of another fractal, respectively (for instance, the Sierpi\'nski carpet or the unit interval, see Figure \ref{fig_VCSC} and Figure \ref{fig_VC01}). The resulting gradient estimates for the heat kernel can behave differently at small and large scales. We will give a detailed discussion of these two examples in Section \ref{sec_eg}.
\end{itemize}
\end{remark}

In this paper, we assume a-priori the strongly recurrent condition (\ref{eq_recurrent}), which singles out specific choices of $\Phi,\Psi$. This framework contains many classical fractals, including the Sierpi\'nski gasket, the Sierpi\'nski carpet, nested fractals, and p.c.f. self-similar sets. Starting from this setting, further examples can be constructed via tensor products. For instance, let $X$ be a space with Hausdorff dimension $\alpha$ and walk dimension $\beta>\alpha$. By taking the $n$-fold tensor product $X^{\otimes n}$, one obtains a space with Hausdorff dimension $n\alpha$ and walk dimension $\beta$. For $n$ sufficiently large, we have $n\alpha\ge\beta$, while the H\"older exponent of harmonic functions on $X^{\otimes n}$ remains equal to $\beta-\alpha$, the same as that on $X$; see \cite[Proposition 3.8]{ABCRST3}. For the general case in which that each factor space is equipped with scaling functions $\Phi,\Psi$, see \cite[Proposition 3.7]{AB25arXiv}.

This paper is organized as follows. In Section \ref{sec_hk}, we list some characterizations of heat kernel estimates for later use. In Section \ref{sec_Poi}, we give the existence, the uniqueness and H\"older estimate for the solutions of Poisson equation. In Section \ref{sec_HR2wBE}, we prove ``\ref{eq_HR}$\Rightarrow$\ref{eq_wBE}". In Section \ref{sec_HR2HHKexp}, we prove ``\ref{eq_HR}$\Rightarrow$\ref{eq_HHKexp}". In Section \ref{sec_wBE2HHK2NLE}, we prove\\
\noindent ``\ref{eq_wBE}$\Rightarrow$\ref{eq_HHK}$\Rightarrow$\ref{eq_NLE}". In Section \ref{sec_NLE2HR}, we prove ``\ref{eq_NLE}$\Rightarrow$\ref{eq_HR}". In Section \ref{sec_eg}, we give some examples related to the Sierpi\'nski carpet and the Vicsek set, including blowups of one fractal by another.

\section{Heat kernel estimates}\label{sec_hk}

In this section, we list some characterizations of heat kernel estimates for later use.

Let $D$ be an open subset of $X$. Denote by $\lambda_1(D)$ the smallest Dirichlet eigenvalue for $D$, that is,
$$\lambda_1(D)=\inf\left\{\frac{\mathcal{E}(u,u)}{\lVert {u}\rVert_2^2}:u\in\mathcal{F}_D\backslash\left\{0\right\}\right\}.$$
We say that the relative Faber-Krahn inequality \ref{eq_FK} holds if there exist $C_F\in(0,+\infty)$, $\nu\in(0,1)$ such that for any ball $B=B(x,r)$, for any open subset $D$ of $B$, we have
\begin{equation*}\label{eq_FK}\tag*{FK($\Psi$)}
\lambda_1(D)\ge\frac{C_F}{\Psi(r)}\left(\frac{m(B)}{m(D)}\right)^{\nu}.
\end{equation*}
We will sometimes use the notation \hypertarget{eq_FKPsinu}{FK($\Psi,\nu$)} to emphasize the role of the value of $\nu$.

We say that the Poincar\'e inequality \ref{eq_PI} holds if there exists $C_P\in(0,+\infty)$ such that for any ball $B=B(x,r)$, for any $u\in\mathcal{F}$, we have
\begin{equation*}\label{eq_PI}\tag*{PI($\Psi$)}
\int_B|u-u_B|^2\mathrm{d}m\le C_{P}\Psi(r)\int_{2B}\mathrm{d}\Gamma(u,u).
\end{equation*}

Let $U,V$ be two open subsets of $X$ satisfying $U\subseteq \overline{{U}}\subseteq V$. We say that $\varphi\in\mathcal{F}$ is a cutoff function for $U\subseteq V$ if $0\le\varphi\le1$ $m$-a.e., $\varphi=1$ $m$-a.e. in an open neighborhood of $\overline{{U}}$ and $\mathrm{supp}(\varphi)\subseteq V$, where $\mathrm{supp}(f)$ refers to the support of the measure of $|f|\mathrm{d}m$ for any given function $f$.

We say that the cutoff Sobolev inequality \ref{eq_CS} holds if there exists $C_S\in(0,+\infty)$ such that for any $x\in X$, for any $R,r\in(0,+\infty)$, there exists a cutoff function $\varphi\in\mathcal{F}$ for $B(x,R)\subseteq B(x,R+r)$ such that for any $f\in\mathcal{F}$, we have
\begin{align*}
&\int_{B(x,R+r)\backslash \overline{{B(x,R)}}}f^2\mathrm{d}\Gamma(\varphi,\varphi)\\
&\le\frac{1}{8}\int_{B(x,R+r)\backslash \overline{{B(x,R)}}}\varphi^2\mathrm{d}\Gamma(f,f)+\frac{C_S}{\Psi(r)}\int_{B(x,R+r)\backslash \overline{{B(x,R)}}}f^2\mathrm{d}m.\label{eq_CS}\tag*{CS($\Psi$)}
\end{align*}

We have the characterizations of heat kernel estimates as follows.

\begin{proposition}[{\cite[Theorem 1.12]{AB15}}]\label{prop_UHK}
Let $(X,d,m,\mathcal{E},\mathcal{F})$ be an unbounded MMD space satisfying \ref{eq_VD}. Then the followings are equivalent.
\begin{enumerate}[label=(\arabic*)]
\item \ref{eq_UHK}.
\item \ref{eq_FK} and \ref{eq_CS}.
\end{enumerate}
\end{proposition}

\begin{proposition}[{\cite[THEOREM 1.2]{GHL15}}]\label{prop_HK}
Let $(X,d,m,\mathcal{E},\mathcal{F})$ be an unbounded MMD space satisfying \ref{eq_VPhi}. Then the followings are equivalent.
\begin{enumerate}[label=(\arabic*)]
\item \ref{eq_UHK} and \ref{eq_NLE}.
\item \ref{eq_PI} and \ref{eq_CS}.
\end{enumerate}
\end{proposition}

\begin{remark}
On any complete non-compact Riemannian manifold, \ref{eq_CS} with $\beta_1=\beta_2=2$, or equivalently, $\Psi(r)=r^2$ for any $r\in(0,+\infty)$, holds automatically, for example, one can simply take
$$\varphi=1-\left(\frac{d(x,\cdot)-(R+\frac{r}{3})}{\frac{r}{3}}\vee0\right)\wedge1,$$
so that the above equivalences hold without \ref{eq_CS} and are classical, see \cite{Gri92,Sal92,Gri94}. However, on a general MMD space, \ref{eq_CS} does not always hold and is involved in the formulation of the previous equivalences.
\end{remark}

\section{Poisson equation on unbounded MMD spaces}\label{sec_Poi}

In this section, we establish existence, uniqueness, and H\"older estimates for solutions of the Poisson equation. Our argument follows a technique analogous to gradient estimate methods for Poisson equations in \cite{DRY23, CJKS20}, see also \cite{Jiang15, JKY14}. The main result of this section is Proposition \ref{prop_Poi_holder}, which will play a central role in the proofs of ``\ref{eq_HR}$\Rightarrow$\ref{eq_wBE}" and ``\ref{eq_HR}$\Rightarrow$\ref{eq_HHKexp}".

First, we have existence, uniqueness, and regularity for solutions of the Poisson equation as follows.

\begin{lemma}[{\cite[Lemma 2.2, Lemma 2.6]{DRY23}}]\label{lem_Poi_exist}
Let $(X,d,m,\mathcal{E},\mathcal{F})$ be an unbounded MMD space satisfying \hyperlink{eq_FKPsinu}{FK($\Psi,\nu$)}. Let $q=\frac{2}{1-\nu}$. Then for any $p\in\left[\frac{q}{q-1},+\infty\right)$, for any ball $B=B(x_0,r)$, for any $f\in L^p(B)$, there exists a unique\footnote{in the sense that if $u_1,u_2\in\mathcal{F}_B$ satisfy $\Delta u_1=\Delta u_2=f$ in $B$, then $u_1=u_2$ $m$-a.e..} $u\in\mathcal{F}_B$ such that $\Delta u=f$ in $B$. There exists $C\in(0,+\infty)$ depending only on $\nu$, $C_F$ such that
$$\dashint_B|u|\mathrm{d}m\le C\sqrt{\Psi(r)}\left(\frac{1}{m(B)}\mathcal{E}(u,u)\right)^{1/2}\le C\Psi(r)\left(\dashint_B|f|^p\mathrm{d}m\right)^{1/p}.$$
\end{lemma}

Second, we have pointwise estimates for solutions of the Poisson equation as follows.

\begin{lemma}[{\cite[Lemma 2.7]{DRY23}}]\label{lem_Poi_point}
Let $(X,d,m,\mathcal{E},\mathcal{F})$ be an unbounded MMD space satisfying \ref{eq_VD}, \hyperlink{eq_FKPsinu}{FK($\Psi,\nu$)} and \ref{eq_CS}. Let $q=\frac{2}{1-\nu}$. Then for any $p\in\left[\frac{q}{q-1},+\infty\right)$, there exists $C\in(0,+\infty)$ such that for any ball $B=B(x_0,r)$, for any $f\in L^\infty(2B)$, if $u\in\mathcal{F}$ satisfies $\Delta u=f$ in $2B$, then for $m$-a.e. $x\in B$, we have
$$|u(x)|\le C\left(\dashint_{2B}|u|\mathrm{d}m+F_1(x)\right),$$
where
$$F_1(x)=\sum_{j\le[\log_2r]}\Psi(2^j)\left(\dashint_{B(x,2^j)}|f|^p\mathrm{d}m\right)^{1/p}.$$ 
\end{lemma}

The proofs of the above two results in \cite{DRY23} were given for the case $\beta_1=2$. It is easy to see that the same proofs remain valid for more general $\Psi$, including the case in our paper.

Third, we have H\"older estimates for solutions of the Poisson equation as follows. The proof is similar to that of \cite[Proposition 5.1]{DRY23} and \cite[Theorem 3.2]{CJKS20} for gradient estimates, with some technical modifications.

\begin{proposition}\label{prop_Poi_holder}
Let $(X,d,m,\mathcal{E},\mathcal{F})$ be an unbounded MMD space satisfying \ref{eq_VPhi},\\
\noindent \hyperlink{eq_FKPsinu}{FK($\Psi,\nu$)}, \ref{eq_CS} and \ref{eq_HR}. Let $q=\frac{2}{1-\nu}$. Then for any $p\in\left[\frac{q}{q-1},+\infty\right)$, there exists $C\in(0,+\infty)$ such that for any ball $B=B(x_0,r)$, for any $f\in L^\infty(2B)$, if $u\in\mathcal{F}$ satisfies $\Delta u=f$ in $2B$, then for $m$-a.e. $x,y\in\frac{1}{16}B$ with $x\ne y$, we have
\begin{align*}
\frac{\Phi(d(x,y))}{\Psi(d(x,y))}|u(x)-u(y)|\le C\left(\frac{\Phi(r)}{\Psi(r)}\dashint_{2B}|u|\mathrm{d}m+F_2(x)+F_2(y)\right),
\end{align*}
where
$$F_2(x)=\sum_{j\le[\log_2r]}\Phi(2^{j})\left(\dashint_{B(x,2^j)}|f|^p\mathrm{d}m\right)^{1/p}.$$
\end{proposition}

\begin{proof}
Let $k_0=[\log_2d(x,y)]$ and $k_1=[\log_2r]$, then $2^{k_0}\le d(x,y)<2^{k_0+1}$ and $2^{k_1}\le r<2^{k_1+1}$. Since $x,y\in\frac{1}{16}B$, we have $d(x,y)<\frac{1}{8}r$. Hence $2^{k_0}\le d(x,y)<\frac{1}{8}r<\frac{1}{8}2^{k_1+1}=2^{k_1-2}$, which implies $k_0+3\le k_1$.

For any $k=k_0+3,\ldots,k_1$, by Lemma \ref{lem_Poi_exist}, there exists a unique $u_k\in\mathcal{F}_{B(x,2^{k})}$ such that $\Delta u_k=f$ in $B(x,2^k)$ and
$$\dashint_{B(x,2^{k-1})}|u_k|\mathrm{d}m\lesssim\dashint_{B(x,2^k)}|u_k|\mathrm{d}m\lesssim\Psi\left(2^{k}\right)\left(\dashint_{B(x,2^{k})}|f|^p\mathrm{d}m\right)^{1/p}.$$
Then
\begin{align*}
&|u(x)-u(y)|\\
&\le|(u-u_{k_1})(x)-(u-u_{k_1})(y)|\\
&+\sum_{k=k_0+4}^{k_1}|(u_k-u_{k-1})(x)-(u_k-u_{k-1})(y)|\\
&+|u_{k_0+3}(x)|+|u_{k_0+3}(y)|.
\end{align*}
For any $k=k_0+3,\ldots,k_1$, we have $d(x,y)<2^{k_0+1}\le2^{k-2}$, that is, $y\in B(x,2^{k-2})$.

Since $\Delta(u-u_{k_1})=0$ in $B(x,2^{k_1})$, by \ref{eq_HR}, we have
\begin{align*}
&\frac{\Phi(d(x,y))}{\Psi(d(x,y))}|(u-u_{k_1})(x)-(u-u_{k_1})(y)|\\
&\lesssim\frac{\Phi(2^{k_1-1})}{\Psi(2^{k_1-1})}\dashint_{B(x,2^{k_1})}|u-u_{k_1}|\mathrm{d}m\\
&\le\frac{\Phi(2^{k_1-1})}{\Psi(2^{k_1-1})}\left(\dashint_{B(x,2^{k_1})}|u|\mathrm{d}m+\dashint_{B(x,2^{k_1})}|u_{k_1}|\mathrm{d}m\right)\\
&\lesssim\frac{\Phi(2^{k_1-1})}{\Psi(2^{k_1-1})}\left(\dashint_{2B}|u|\mathrm{d}m+\Psi\left(2^{k_1}\right)\left(\dashint_{B(x,2^{k_1})}|f|^p\mathrm{d}m\right)^{1/p}\right)\\
&\lesssim\frac{\Phi(r)}{\Psi(r)}\dashint_{2B}|u|\mathrm{d}m+\Phi(2^{k_1})\left(\dashint_{B(x,2^{k_1})}|f|^p\mathrm{d}m\right)^{1/p}.
\end{align*}

Similarly, for any $k=k_0+4,\ldots,k_1$, since $\Delta(u_k-u_{k-1})=0$ in $B(x,2^{k-1})$, by \ref{eq_HR}, we have
\begin{align*}
&\frac{\Phi(d(x,y))}{\Psi(d(x,y))}|(u_k-u_{k-1})(x)-(u_k-u_{k-1})(y)|\\
&\lesssim\frac{\Phi(2^{k-2})}{\Psi(2^{k-2})}\dashint_{B(x,2^{k-1})}|u_k-u_{k-1}|\mathrm{d}m\\
&\le\frac{\Phi(2^{k-2})}{\Psi(2^{k-2})}\left(\dashint_{B(x,2^{k-1})}|u_k|\mathrm{d}m+\dashint_{B(x,2^{k-1})}|u_{k-1}|\mathrm{d}m\right)\\
&\lesssim\frac{\Phi(2^{k-2})}{\Psi(2^{k-2})}\left(\Psi\left(2^{k}\right)\left(\dashint_{B(x,2^k)}|f|^p\mathrm{d}m\right)^{1/p}+\Psi\left(2^{k-1}\right)\left(\dashint_{B(x,2^{k-1})}|f|^p\mathrm{d}m\right)^{1/p}\right)\\
&\lesssim\Phi(2^{k})\left(\dashint_{B(x,2^{k})}|f|^p\mathrm{d}m\right)^{1/p}+\Phi(2^{k-1})\left(\dashint_{B(x,2^{k-1})}|f|^p\mathrm{d}m\right)^{1/p}.
\end{align*}

Since $\Delta u_{k_0+3}=f$ in $B(x,2^{k_0+3})$, by Lemma \ref{lem_Poi_point}, we have
\begin{align*}
&\frac{\Phi(d(x,y))}{\Psi(d(x,y))}|u_{k_0+3}(x)|\\
&\lesssim\frac{\Phi(d(x,y))}{\Psi(d(x,y))}\left(\dashint_{B(x,2^{k_0+3})}|u_{k_0+3}|\mathrm{d}m+\sum_{j\le k_0+2}\Psi(2^j)\left(\dashint_{B(x,2^{j})}|f|^p\mathrm{d}m\right)^{1/p}\right)\\
&\lesssim\frac{\Phi(d(x,y))}{\Psi(d(x,y))}\left(\Psi(2^{k_0+3})\left(\dashint_{B(x,2^{k_0+3})}|f|^p\mathrm{d}m\right)^{1/p}+\sum_{j\le k_0+2}\Psi(2^j)\left(\dashint_{B(x,2^{j})}|f|^p\mathrm{d}m\right)^{1/p}\right)\\
&\lesssim\frac{\Phi(2^{k_0+3})}{\Psi(2^{k_0+3})}\Psi(2^{k_0+3})\left(\dashint_{B(x,2^{k_0+3})}|f|^p\mathrm{d}m\right)^{1/p}\\
&+\sum_{j\le k_0+2}\frac{\Phi(2^{k_0+2})}{\Psi(2^{k_0+2})}\frac{\Psi(2^j)}{\Phi(2^j)}\Phi(2^j)\left(\dashint_{B(x,2^{j})}|f|^p\mathrm{d}m\right)^{1/p}\\
&\le\Phi(2^{k_0+3})\left(\dashint_{B(x,2^{k_0+3})}|f|^p\mathrm{d}m\right)^{1/p}+\sum_{j\le k_0+2}\Phi(2^j)\left(\dashint_{B(x,2^{j})}|f|^p\mathrm{d}m\right)^{1/p},
\end{align*}
where in the last inequality, we use the fact that $\gamma_1,\gamma_2>0$, which implies that
$$\frac{\Phi(2^{k_0+2})}{\Psi(2^{k_0+2})}\frac{\Psi(2^j)}{\Phi(2^j)}\le1\text{ for any }j\le k_0+2.$$
Similarly, since $y\in B(x,2^{k_0+1})\subseteq B(x,2^{k_0+2})$, by Lemma \ref{lem_Poi_point}, we also have
\begin{align*}
&\frac{\Phi(d(x,y))}{\Psi(d(x,y))}|u_{k_0+3}(y)|\\
&\lesssim\frac{\Phi(d(x,y))}{\Psi(d(x,y))}\left(\dashint_{B(x,2^{k_0+3})}|u_{k_0+3}|\mathrm{d}m+\sum_{j\le k_0+2}\Psi(2^j)\left(\dashint_{B(y,2^{j})}|f|^p\mathrm{d}m\right)^{1/p}\right)\\
&\lesssim\Phi(2^{k_0+3})\left(\dashint_{B(x,2^{k_0+3})}|f|^p\mathrm{d}m\right)^{1/p}+\sum_{j\le k_0+2}\Phi(2^j)\left(\dashint_{B(y,2^{j})}|f|^p\mathrm{d}m\right)^{1/p}.
\end{align*}

In summary, we have
\begin{align*}
&\frac{\Phi(d(x,y))}{\Psi(d(x,y))}|u(x)-u(y)|\\
&\lesssim\frac{\Phi(r)}{\Psi(r)}\dashint_{2B}|u|\mathrm{d}m+\Phi(2^{k_1})\left(\dashint_{B(x,2^{k_1})}|f|^p\mathrm{d}m\right)^{1/p}\\
&+\sum_{k=k_0+4}^{k_1}\left(\Phi(2^{k})\left(\dashint_{B(x,2^{k})}|f|^p\mathrm{d}m\right)^{1/p}+\Phi(2^{k-1})\left(\dashint_{B(x,2^{k-1})}|f|^p\mathrm{d}m\right)^{1/p}\right)\\
&+\Phi(2^{k_0+3})\left(\dashint_{B(x,2^{k_0+3})}|f|^p\mathrm{d}m\right)^{1/p}+\sum_{j\le k_0+2}\Phi(2^j)\left(\dashint_{B(x,2^{j})}|f|^p\mathrm{d}m\right)^{1/p}\\
&+\Phi(2^{k_0+3})\left(\dashint_{B(x,2^{k_0+3})}|f|^p\mathrm{d}m\right)^{1/p}+\sum_{j\le k_0+2}\Phi(2^j)\left(\dashint_{B(y,2^{j})}|f|^p\mathrm{d}m\right)^{1/p}\\
&\lesssim\frac{\Phi(r)}{\Psi(r)}\dashint_{2B}|u|\mathrm{d}m+\sum_{k=k_0+3}^{k_1}\Phi(2^{k})\left(\dashint_{B(x,2^{k})}|f|^p\mathrm{d}m\right)^{1/p}\\
&+\sum_{k\le k_0+3}\Phi(2^k)\left(\dashint_{B(x,2^{k})}|f|^p\mathrm{d}m\right)^{1/p}+\sum_{k\le k_0+2}\Phi(2^k)\left(\dashint_{B(y,2^{k})}|f|^p\mathrm{d}m\right)^{1/p}\\
&\lesssim\frac{\Phi(r)}{\Psi(r)}\dashint_{2B}|u|\mathrm{d}m+\sum_{k\le k_1}\Phi(2^{k})\left(\dashint_{B(x,2^{k})}|f|^p\mathrm{d}m\right)^{1/p}+\sum_{k\le k_1}\Phi(2^k)\left(\dashint_{B(y,2^{k})}|f|^p\mathrm{d}m\right)^{1/p}.
\end{align*}
\end{proof}

\section{Proof of \texorpdfstring{``\ref{eq_HR}$\Rightarrow$\ref{eq_wBE}"}{"HR implies wBE"}}\label{sec_HR2wBE}

The idea is to consider the heat equation ``$\Delta(P_tf)=\frac{\partial}{\partial t}(P_tf)$" for fixed $t$ as a Poisson equation.

\begin{proof}[Proof of ``\ref{eq_HR}$\Rightarrow$\ref{eq_wBE}"]
For any $t\in(0,+\infty)$, for any $f\in L^\infty(X;m)$, we have $P_tf$ is a solution of the Poisson equation $\Delta(P_tf)=\frac{\partial}{\partial t}(P_tf)$. For any $x,y\in X$ with $x\ne y$. If $d(x,y)<\Psi^{-1}(t)$, then taking $r=16\Psi^{-1}(t)>16d(x,y)$, by Proposition \ref{prop_Poi_holder}, we have
\begin{align*}
&\frac{\Phi(d(x,y))}{\Psi(d(x,y))}|P_tf(x)-P_tf(y)|\\
&\lesssim\frac{\Phi(16\Psi^{-1}(t))}{\Psi(16\Psi^{-1}(t))}\dashint_{B(x,32\Psi^{-1}(t))}|P_tf|\mathrm{d}m\\
&+\sum_{k\le[\log_216\Psi^{-1}(t)]}\Phi(2^k)\left(\dashint_{B(x,2^k)}|\frac{\partial}{\partial t}(P_tf)|^p\mathrm{d}m\right)^{1/p}\\
&+\sum_{k\le[\log_216\Psi^{-1}(t)]}\Phi(2^k)\left(\dashint_{B(y,2^k)}|\frac{\partial}{\partial t}(P_tf)|^p\mathrm{d}m\right)^{1/p}.
\end{align*}
For the first term, we have
\begin{align*}
&\frac{\Phi(16\Psi^{-1}(t))}{\Psi(16\Psi^{-1}(t))}\dashint_{B(x,32\Psi^{-1}(t))}|P_tf|\mathrm{d}m\\
&\lesssim\frac{\Phi(\Psi^{-1}(t))}{\Psi(\Psi^{-1}(t))}\dashint_{B(x,32\Psi^{-1}(t))}\lVert {P_tf}\rVert_{L^\infty(X;m)}\mathrm{d}m\\
&\le\frac{\Phi(\Psi^{-1}(t))}{\Psi(\Psi^{-1}(t))}\lVert {f}\rVert_{L^\infty(X;m)}.
\end{align*}
For the second term, by \ref{eq_UHK} and \cite[THEOREM 4]{Dav97}, we have the following estimate of the time derivative of the heat kernel
\begin{equation*}
\left|\frac{\partial}{\partial t}p_t(z,w)\right|\le\frac{C_1}{tV\left(z,\Psi^{-1}(t)\right)}\exp\left(-\Upsilon\left(C_2d(z,w),t\right)\right),
\end{equation*}
hence
\begin{align*}
&|\frac{\partial}{\partial t}(P_tf)(z)|\\
&\le\int_X\lvert\frac{\partial}{\partial t}p_t(z,w)\rvert\cdot|f(w)|m(\mathrm{d}w)\\
&\le \lVert {f}\rVert_{L^\infty(X;m)}\int_X\frac{C_1}{tV\left(z,\Psi^{-1}(t)\right)}\exp\left(-\Upsilon\left(C_2d(z,w),t\right)\right)m(\mathrm{d}w).
\end{align*}
By writing
$$\int_X=\int_{B(z,\Psi^{-1}(t))}+\sum_{n=0}^\infty\int_{B(z,2^{n+1}\Psi^{-1}(t))\backslash B(z,2^n\Psi^{-1}(t))},$$
and using \ref{eq_VPhi}, we have the integral
$$\int_X\frac{1}{V\left(z,\Psi^{-1}(t)\right)}\exp\left(-\Upsilon\left(C_2d(z,w),t\right)\right)m(\mathrm{d}w)$$
is bounded by some constant depending only on $C_2$ and $\Phi$, see also \cite[Page 944, Line -5]{ACDH04}, \cite[Equation (2.5)]{CCFR17}, hence
$$|\frac{\partial}{\partial t}(P_tf)(z)|\lesssim\frac{\lVert {f}\rVert_{L^\infty(X;m)}}{t}.$$
Hence
\begin{align*}
&\sum_{k\le[\log_216\Psi^{-1}(t)]}\Phi(2^k)\left(\dashint_{B(x,2^k)}|\frac{\partial}{\partial t}(P_tf)|^p\mathrm{d}m\right)^{1/p}\\
&\lesssim\sum_{k\le[\log_216\Psi^{-1}(t)]}\Phi(2^k)\left(\dashint_{B(x,2^k)}\left(\frac{\lVert {f}\rVert_{L^\infty(X;m)}}{t}\right)^p\mathrm{d}m\right)^{1/p}\\
&=\frac{\lVert{f} \rVert_{L^\infty(X;m)}}{t}\sum_{k\le[\log_216\Psi^{-1}(t)]}\frac{\Phi(2^k)}{\Phi(16\Psi^{-1}(t))}\frac{\Phi(16\Psi^{-1}(t))}{\Phi(\Psi^{-1}(t))}\Phi(\Psi^{-1}(t))\\
&\lesssim\frac{\lVert {f}\rVert_{L^\infty(X;m)}}{t}\Phi(\Psi^{-1}(t))\\
&=\frac{\Phi(\Psi^{-1}(t))}{\Psi(\Psi^{-1}(t))}\lVert {f}\rVert_{L^\infty(X;m)}.
\end{align*}
For the third term, similarly, we have
$$\sum_{k\le[\log_216\Psi^{-1}(t)]}\Phi(2^k)\left(\dashint_{B(y,2^k)}|\frac{\partial}{\partial t}(P_tf)|^p\mathrm{d}m\right)^{1/p}\lesssim\frac{\Phi(\Psi^{-1}(t))}{\Psi(\Psi^{-1}(t))}\lVert{f} \rVert_{L^\infty(X;m)}.$$
In summary, we have
$$\frac{\Phi(d(x,y))}{\Psi(d(x,y))}|P_tf(x)-P_tf(y)|\lesssim\frac{\Phi(\Psi^{-1}(t))}{\Psi(\Psi^{-1}(t))}\lVert {f}\rVert_{L^\infty(X;m)}.$$

If $d(x,y)\ge\Psi^{-1}(t)$, then obviously
\begin{align*}
&\frac{\Phi(d(x,y))}{\Psi(d(x,y))}|P_tf(x)-P_tf(y)|\\
&\le2\frac{\Phi(d(x,y))}{\Psi(d(x,y))}\lVert {P_tf}\rVert_{L^\infty(X;m)}\le2\frac{\Phi(d(x,y))}{\Psi(d(x,y))}\lVert {f}\rVert_{L^\infty(X;m)}\\
&=2\frac{\frac{\Psi(\Psi^{-1}(t))}{\Phi(\Psi^{-1}(t))}}{\frac{\Psi(d(x,y))}{\Phi(d(x,y))}}\frac{\Phi(\Psi^{-1}(t))}{\Psi(\Psi^{-1}(t))}\lVert {f}\rVert_{L^\infty(X;m)}\le2\frac{\Phi(\Psi^{-1}(t))}{\Psi(\Psi^{-1}(t))}\lVert {f}\rVert_{L^\infty(X;m)}.
\end{align*}
\end{proof}

\section{Proof of \texorpdfstring{``\ref{eq_HR}$\Rightarrow$\ref{eq_HHKexp}"}{"HR implies HHKexp"}}\label{sec_HR2HHKexp}

The idea is to consider the heat equation ``$\Delta p_t(x,\cdot)=\frac{\partial}{\partial t}p_t(x,\cdot)$" for fixed $t,x$ as a Poisson equation. Comparing with the proof of ``\ref{eq_HR}$\Rightarrow$\ref{eq_wBE}" in Section \ref{sec_HR2wBE}, we need to pay special attention to the exponential terms.

\begin{proof}[Proof of ``\ref{eq_HR}$\Rightarrow$\ref{eq_HHKexp}"]
By \ref{eq_UHK} and \cite[THEOREM 4]{Dav97}, we have the following estimates of the heat kernel
$$p_t(x,y)\le\frac{C_1}{V(x,\Psi^{-1}(t))}\exp\left(-\Upsilon(C_2d(x,y),t)\right),$$
and its time derivative
$$\left|\frac{\partial}{\partial t}p_t(z,w)\right|\le\frac{C_1}{tV\left(z,\Psi^{-1}(t)\right)}\exp\left(-\Upsilon\left(C_2d(z,w),t\right)\right).$$

Take any $t\in(0,+\infty)$, $x,y_1,y_2\in X$. We now distinguish the following four cases:
\begin{enumerate}[label=(\roman*),ref=(\roman*)]
\item\label{HHKcase1} $d(y_1,y_2)\ge\Psi^{-1}(t)$.
\item\label{HHKcase2} $d(y_1,y_2)<\Psi^{-1}(t)$, $d(x,y_1)<128\Psi^{-1}(t)$ and $d(x,y_2)<128\Psi^{-1}(t)$.
\item\label{HHKcase3} $d(y_1,y_2)<\Psi^{-1}(t)$ and $d(x,y_1)\ge128\Psi^{-1}(t)$.
\item\label{HHKcase4} $d(y_1,y_2)<\Psi^{-1}(t)$ and $d(x,y_2)\ge128\Psi^{-1}(t)$.
\end{enumerate}

For case \ref{HHKcase1}, we obtain the result directly from \ref{eq_UHK} as follows.
\begin{align*}
&|p_t(x,y_1)-p_t(x,y_2)|\le p_t(x,y_1)+p_t(x,y_2)\\
&\le\frac{C_1}{V(x,\Psi^{-1}(t))}\exp\left(-\Upsilon(C_2d(x,y_1),t)\right)+\frac{C_1}{V(x,\Psi^{-1}(t))}\exp\left(-\Upsilon(C_2d(x,y_2),t)\right)\\
&\asymp\frac{1}{\Phi(\Psi^{-1}(t))}\left(\exp\left(-\Upsilon(C_2d(x,y_1),t)\right)+\exp\left(-\Upsilon(C_2d(x,y_2),t)\right)\right)\\
&=\frac{\Psi(\Psi^{-1}(t))}{\Phi(\Psi^{-1}(t))}\frac{1}{t}\left(\exp\left(-\Upsilon(C_2d(x,y_1),t)\right)+\exp\left(-\Upsilon(C_2d(x,y_2),t)\right)\right)\\
&\le\frac{\left(\frac{\Psi}{\Phi}\right)(d(y_1,y_2))}{t}\left(\exp\left(-\Upsilon(C_2d(x,y_1),t)\right)+\exp\left(-\Upsilon(C_2d(x,y_2),t)\right)\right).
\end{align*}

For the remaining three cases, by Proposition \ref{prop_Poi_holder}, letting $r=16\Psi^{-1}(t)>16d(y_1,y_2)$, we have
\begin{align}
&\frac{\Phi(d(y_1,y_2))}{\Psi(d(y_1,y_2))}|p_t(x,y_1)-p_t(x,y_2)|\nonumber\\
&\lesssim\frac{\Phi(16\Psi^{-1}(t))}{\Psi(16\Psi^{-1}(t))}\dashint_{B(y_1,32\Psi^{-1}(t))}p_t(x,z)m(\mathrm{d}z)\nonumber\\
&+\sum_{j\le[\log_216\Psi^{-1}(t)]}\Phi(2^j)\left(\dashint_{B(y_1,2^j)}\left|\frac{\partial}{\partial t}p_t(x,z)\right|^pm(\mathrm{d}z)\right)^{1/p}\nonumber\\
&+\sum_{j\le[\log_216\Psi^{-1}(t)]}\Phi(2^j)\left(\dashint_{B(y_2,2^j)}\left|\frac{\partial}{\partial t}p_t(x,z)\right|^pm(\mathrm{d}z)\right)^{1/p}.\label{eq_HHKexp0}
\end{align}

For case \ref{HHKcase2}, we have $d(x,y_1)<128\Psi^{-1}(t)$ and $d(x,y_2)<128\Psi^{-1}(t)$. By Equation (\ref{eq_HHKexp0}), we have
\begin{align*}
&\frac{\Phi(d(y_1,y_2))}{\Psi(d(y_1,y_2))}|p_t(x,y_1)-p_t(x,y_2)|\\
&\lesssim\frac{\Phi(16\Psi^{-1}(t))}{\Psi(16\Psi^{-1}(t))}\dashint_{B(y_1,32\Psi^{-1}(t))}\frac{C_1}{V(x,\Psi^{-1}(t))}m(\mathrm{d}z)\\
&+\sum_{j\le[\log_216\Psi^{-1}(t)]}\Phi(2^j)\left(\dashint_{B(y_1,2^j)}\left(\frac{C_1}{tV(x,\Psi^{-1}(t))}\right)^pm(\mathrm{d}z)\right)^{1/p}\\
&+\sum_{j\le[\log_216\Psi^{-1}(t)]}\Phi(2^j)\left(\dashint_{B(y_2,2^j)}\left(\frac{C_1}{tV(x,\Psi^{-1}(t))}\right)^pm(\mathrm{d}z)\right)^{1/p}\\
&=\frac{\Phi(16\Psi^{-1}(t))}{\Psi(16\Psi^{-1}(t))}\frac{C_1}{V(x,\Psi^{-1}(t))}+2\sum_{j\le[\log_216\Psi^{-1}(t)]}\Phi(2^j)\frac{C_1}{tV(x,\Psi^{-1}(t))}\\
&\asymp\frac{\Phi(\Psi^{-1}(t))}{\Psi(\Psi^{-1}(t))}\frac{1}{\Phi(\Psi^{-1}(t))}+\frac{1}{t\Phi(\Psi^{-1}(t))}\sum_{j\le[\log_216\Psi^{-1}(t)]}\Phi(2^j)\\
&=\frac{1}{t}+\frac{1}{t\Phi(\Psi^{-1}(t))}\sum_{j\le[\log_216\Psi^{-1}(t)]}\Phi(2^j).
\end{align*}

For any $r\in(0,+\infty)$. If $r<1$, then
$$\sum_{j\le[\log_2r]}\Phi(2^j)=\sum_{j\le[\log_2r]}(2^j)^{\alpha_1}\asymp2^{[\log_2r]\alpha_1}\asymp r^{\alpha_1}=\Phi(r).$$
If $r\ge1$, then
$$\sum_{j\le[\log_2r]}\Phi(2^j)=\sum_{j\le-1}(2^j)^{\alpha_1}+\sum_{j=0}^{[\log_2r]}(2^j)^{\alpha_2}\asymp1+2^{[\log_2r]\alpha_2}\asymp1+r^{\alpha_2}\asymp r^{\alpha_2}=\Phi(r).$$
Hence
\begin{equation*}
\sum_{j\le[\log_2r]}\Phi(2^j)\asymp\Phi(r)\text{ for any }r\in(0,+\infty).
\end{equation*}

Then
\begin{align*}
&\frac{\Phi(d(y_1,y_2))}{\Psi(d(y_1,y_2))}|p_t(x,y_1)-p_t(x,y_2)|\lesssim\frac{1}{t}+\frac{1}{t\Phi(\Psi^{-1}(t))}\Phi(16\Psi^{-1}(t))\\
&\asymp\frac{1}{t}+\frac{1}{t\Phi(\Psi^{-1}(t))}\Phi(\Psi^{-1}(t))=\frac{2}{t}.
\end{align*}
By assumption, for $i=1,2$, we have
$$\Upsilon(C_2d(x,y_i),t)=\sup_{s\in(0,+\infty)}\left(\frac{C_2d(x,y_i)}{\Psi^{-1}(s)}-\frac{t}{s}\right)\le\sup_{s\in(0,+\infty)}\left(\frac{128C_2\Psi^{-1}(t)}{\Psi^{-1}(s)}-\frac{t}{s}\right).$$
By the following elementary result, Lemma \ref{lem_HHKexp_ele}, there exists $C_3\in(0,+\infty)$ depending only on $C_2$, $\beta_1$, $\beta_2$ such that
$$\sup_{t,s\in(0,+\infty)}\left(\frac{128C_2\Psi^{-1}(t)}{\Psi^{-1}(s)}-\frac{t}{s}\right)\le C_3.$$
Hence
\begin{align*}
&\frac{\Phi(d(y_1,y_2))}{\Psi(d(y_1,y_2))}|p_t(x,y_1)-p_t(x,y_2)|\lesssim\frac{e^{C_3}}{t}(2e^{-C_3})\\
&\le\frac{e^{C_3}}{t}\left(\exp(-\Upsilon(C_2d(x,y_1),t))+\exp(-\Upsilon(C_2d(x,y_2),t))\right).
\end{align*}

For case \ref{HHKcase3}, we have $d(x,y_1)\ge128\Psi^{-1}(t)$. Then in Equation (\ref{eq_HHKexp0}), for any $z\in B(y_1,32\Psi^{-1}(t))$, we have
\begin{align*}
&d(x,z)\ge d(x,y_1)-d(y_1,z)\ge d(x,y_1)-32\Psi^{-1}(t)\\
&\ge d(x,y_1)-\frac{1}{4}d(x,y_1)=\frac{3}{4}d(x,y_1)\ge\frac{1}{2}d(x,y_1).
\end{align*}
For any $j\le[\log_216\Psi^{-1}(t)]$, for any $z\in B(y_1,2^j)$, we have
\begin{align*}
&d(x,z)\ge d(x,y_1)-d(y_1,z)\ge d(x,y_1)-16\Psi^{-1}(t)\\
&\ge d(x,y_1)-\frac{1}{8}d(x,y_1)=\frac{7}{8}d(x,y_1)\ge\frac{1}{2}d(x,y_1),
\end{align*}
for any $z\in B(y_2,2^j)$, we have
\begin{align*}
&d(x,z)\ge d(x,y_1)-d(y_1,y_2)-d(y_2,z)\ge d(x,y_1)-\Psi^{-1}(t)-16\Psi^{-1}(t)\\
&\ge d(x,y_1)-\frac{1}{4}d(x,y_1)=\frac{3}{4}d(x,y_1)\ge\frac{1}{2}d(x,y_1).
\end{align*}
Hence Equation (\ref{eq_HHKexp0}) gives
\begin{align*}
&\frac{\Phi(d(y_1,y_2))}{\Psi(d(y_1,y_2))}|p_t(x,y_1)-p_t(x,y_2)|\\
&\lesssim\frac{\Phi(16\Psi^{-1}(t))}{\Psi(16\Psi^{-1}(t))}\dashint_{B(y_1,32\Psi^{-1}(t))}\frac{C_1}{V(x,\Psi^{-1}(t))}\exp\left(-\Upsilon(C_2d(x,z),t)\right)m(\mathrm{d}z)\\
&+\sum_{j\le[\log_216\Psi^{-1}(t)]}\Phi(2^j)\left(\dashint_{B(y_1,2^j)}\left(\frac{C_1}{tV(x,\Psi^{-1}(t))}\exp\left(-\Upsilon(C_2d(x,z),t)\right)\right)^pm(\mathrm{d}z)\right)^{1/p}\\
&+\sum_{j\le[\log_216\Psi^{-1}(t)]}\Phi(2^j)\left(\dashint_{B(y_2,2^j)}\left(\frac{C_1}{tV(x,\Psi^{-1}(t))}\exp\left(-\Upsilon(C_2d(x,z),t)\right)\right)^pm(\mathrm{d}z)\right)^{1/p}\\
&\le\frac{\Phi(16\Psi^{-1}(t))}{\Psi(16\Psi^{-1}(t))}\dashint_{B(y_1,32\Psi^{-1}(t))}\frac{C_1}{V(x,\Psi^{-1}(t))}\exp\left(-\Upsilon(\frac{C_2}{2}d(x,y_1),t)\right)m(\mathrm{d}z)\\
&+\sum_{j\le[\log_216\Psi^{-1}(t)]}\Phi(2^j)\left(\dashint_{B(y_1,2^j)}\left(\frac{C_1}{tV(x,\Psi^{-1}(t))}\exp\left(-\Upsilon(\frac{C_2}{2}d(x,y_1),t)\right)\right)^pm(\mathrm{d}z)\right)^{1/p}\\
&+\sum_{j\le[\log_216\Psi^{-1}(t)]}\Phi(2^j)\left(\dashint_{B(y_2,2^j)}\left(\frac{C_1}{tV(x,\Psi^{-1}(t))}\exp\left(-\Upsilon(\frac{C_2}{2}d(x,y_1),t)\right)\right)^pm(\mathrm{d}z)\right)^{1/p}\\
&=\frac{\Phi(16\Psi^{-1}(t))}{\Psi(16\Psi^{-1}(t))}\frac{C_1}{V(x,\Psi^{-1}(t))}\exp\left(-\Upsilon(\frac{C_2}{2}d(x,y_1),t)\right)\\
&+2\sum_{j\le[\log_216\Psi^{-1}(t)]}\Phi(2^j)\frac{C_1}{tV(x,\Psi^{-1}(t))}\exp\left(-\Upsilon(\frac{C_2}{2}d(x,y_1),t)\right)\\
&\asymp\frac{\Phi(\Psi^{-1}(t))}{\Psi(\Psi^{-1}(t))}\frac{1}{\Phi(\Psi^{-1}(t))}\exp\left(-\Upsilon(\frac{C_2}{2}d(x,y_1),t)\right)\\
&+\frac{1}{t\Phi(\Psi^{-1}(t))}\exp\left(-\Upsilon(\frac{C_2}{2}d(x,y_1),t)\right)\sum_{j\le[\log_216\Psi^{-1}(t)]}\Phi(2^j)\\
&=\frac{1}{t}\left(1+\frac{1}{\Phi(\Psi^{-1}(t))}\sum_{j\le[\log_216\Psi^{-1}(t)]}\Phi(2^j)\right)\exp\left(-\Upsilon(\frac{C_2}{2}d(x,y_1),t)\right)\\
&\asymp\frac{1}{t}\left(1+\frac{1}{\Phi(\Psi^{-1}(t))}\Phi(16\Psi^{-1}(t))\right)\exp\left(-\Upsilon(\frac{C_2}{2}d(x,y_1),t)\right)\\
&\asymp\frac{1}{t}\exp\left(-\Upsilon(\frac{C_2}{2}d(x,y_1),t)\right).
\end{align*}

For case \ref{HHKcase4}, similar to case \ref{HHKcase3}, we also have
$$\frac{\Phi(d(y_1,y_2))}{\Psi(d(y_1,y_2))}|p_t(x,y_1)-p_t(x,y_2)|\lesssim\frac{1}{t}\exp\left(-\Upsilon(\frac{C_2}{2}d(x,y_2),t)\right).$$

In summary, for any $t\in(0,+\infty)$, $x,y_1,y_2\in X$, we have
\begin{align*}
&|p_t(x,y_1)-p_t(x,y_2)|\\
&\lesssim\frac{1}{t}\frac{\Psi(d(y_1,y_2))}{\Phi(d(y_1,y_2))}\left(\exp\left(-\Upsilon(\frac{C_2}{2}d(x,y_1),t)\right)+\exp\left(-\Upsilon(\frac{C_2}{2}d(x,y_2),t)\right)\right).
\end{align*}
\end{proof}

\begin{lemma}[{\cite[Lemma 5.2]{DRY23}}]\label{lem_HHKexp_ele}
Let $A\in(0,+\infty)$. Then there exists $C\in(0,+\infty)$ depending only on $A$, $\beta_1$, $\beta_2$ such that
$$\sup_{t,s\in(0,+\infty)}\left(A\frac{\Psi^{-1}(t)}{\Psi^{-1}(s)}-\frac{t}{s}\right)\le C.$$
\end{lemma}

\begin{proof}
Let $f:(0,+\infty)\times(0,+\infty)\to \mathbb{R}$ be given by
$$f(t,s)=A\frac{\Psi^{-1}(t)}{\Psi^{-1}(s)}-\frac{t}{s}.$$
By considering separately the cases $t,s\in(0,1)$, $t,s\in[1,+\infty)$, $0<t<1\le s$ and $0<s<1\le t$, we easily obtain that
$$f(t,s)\le A\max\left\{\left(\frac{t}{s}\right)^{1/\beta_1},\left(\frac{t}{s}\right)^{1/\beta_2}\right\}-\frac{t}{s}.$$
Since the function $(0,+\infty)\to \mathbb{R}$, $x\mapsto A\max\{x^{1/\beta_1},x^{1/\beta_2}\}-x$ is bounded from above by some positive constant $C$ depending only on $A$, $\beta_1$, $\beta_2$, we have
$$\sup_{t,s\in(0,+\infty)}f(t,s)\le C.$$
\end{proof}

\section{Proof of \texorpdfstring{``\ref{eq_wBE}$\Rightarrow$\ref{eq_HHK}$\Rightarrow$\ref{eq_NLE}"}{"wBE implies HHK implies NLE"}}\label{sec_wBE2HHK2NLE}

The proof of ``\ref{eq_wBE}$\Rightarrow$\ref{eq_HHK}" is direct by using the semi-group property, see also \cite[Lemma 3.4]{ABCRST3}.

\begin{proof}[Proof of ``\ref{eq_wBE}$\Rightarrow$\ref{eq_HHK}"]
For any $t\in(0,+\infty)$, for any $x,y_1,y_2\in X$, we have
\begin{align*}
&p_t(x,y_1)=p_t(y_1,x)=\int_Xp_{t/2}(y_1,z)p_{t/2}(z,x)m(\mathrm{d}z)\\
&=\int_Xp_{t/2}(y_1,z)p_{t/2}(x,z)m(\mathrm{d}z)=P_{t/2}(p_{t/2}(x,\cdot))(y_1),
\end{align*}
and $p_t(x,y_2)=P_{t/2}(p_{t/2}(x,\cdot))(y_2)$. By \ref{eq_wBE} and \ref{eq_UHK}, we have
\begin{align*}
&\frac{\Phi(d(y_1,y_2))}{\Psi(d(y_1,y_2))}|p_t(x,y_1)-p_t(x,y_2)|\\
&=\frac{\Phi(d(y_1,y_2))}{\Psi(d(y_1,y_2))}|P_{t/2}(p_{t/2}(x,\cdot))(y_1)-P_{t/2}(p_{t/2}(x,\cdot))(y_2)|\\
&\le C_{wBE}\frac{\Phi(\Psi^{-1}(\frac{t}{2}))}{\Psi(\Psi^{-1}(\frac{t}{2}))}\lVert {p_{t/2}(x,\cdot)}\rVert_{L^\infty(X;m)}\\
&\le C_{wBE}\frac{\Phi(\Psi^{-1}(\frac{t}{2}))}{\frac{t}{2}}\frac{C_1}{V(x,\Psi^{-1}(t/2))}\\
&\le 2C_{wBE}\frac{\Phi(\Psi^{-1}(\frac{t}{2}))}{t}\frac{C_1C_{VR}}{\Phi(\Psi^{-1}(t/2))}\\
&=\left(2C_1C_{VR}C_{wBE}\right)\frac{1}{t},
\end{align*}
that is,
$$|p_t(x,y_1)-p_t(x,y_2)|\lesssim\frac{1}{t}\frac{\Psi(d(y_1,y_2))}{\Phi(d(y_1,y_2))}.$$
Hence for any $x_1,x_2\in X$, we have
\begin{align*}
&|p_t(x_1,y_1)-p_t(x_2,y_2)|\\
&\le|p_t(x_1,y_1)-p_t(x_1,y_2)|+|p_t(x_1,y_2)-p_t(x_2,y_2)|\\
&\lesssim\frac{\left(\frac{\Psi}{\Phi}\right)(d(y_1,y_2))+\left(\frac{\Psi}{\Phi}\right)(d(x_1,x_2))}{t}.
\end{align*}
\end{proof}

The proof of ``\ref{eq_HHK}$\Rightarrow$\ref{eq_NLE}" is classical and standard, see also \cite[Theorem 3.1]{Cou03}, \cite[Section 5]{GT12}.

\begin{proof}[Proof of ``\ref{eq_HHK}$\Rightarrow$\ref{eq_NLE}"]
First, we give a lower bound of $p_t(x,x)$. For any $N\ge0$, by \ref{eq_UHK}, we have
\begin{align*}
&\int_{X\backslash B(x,2^N\Psi^{-1}(t))}p_t(x,y)m(\mathrm{d}y)\\
&=\sum_{n=N}^\infty\int_{B(x,2^{n+1}\Psi^{-1}(t))\backslash B(x,2^n\Psi^{-1}(t))}p_t(x,y)m(\mathrm{d}y)\\
&\le\sum_{n=N}^\infty\int_{B(x,2^{n+1}\Psi^{-1}(t))\backslash B(x,2^n\Psi^{-1}(t))}\frac{C_1}{V(x,\Psi^{-1}(t))}\exp\left(-\Upsilon(C_2d(x,y),t)\right)m(\mathrm{d}y)\\
&\le\sum_{n=N}^\infty C_1\frac{V(x,2^{n+1}\Psi^{-1}(t))}{V(x,\Psi^{-1}(t))}\exp\left(-\Upsilon(C_22^n\Psi^{-1}(t),t)\right).
\end{align*}
By \ref{eq_VPhi}, we have
$$\frac{V(x,2^{n+1}\Psi^{-1}(t))}{V(x,\Psi^{-1}(t))}\le C_{VR}^2(2^{n+1})^{\alpha_1\vee\alpha_2}.$$
Note that
\begin{align*}
&\Upsilon(C_22^n\Psi^{-1}(t),t)=\sup_{s\in(0,+\infty)}\left(\frac{C_22^n\Psi^{-1}(t)}{s}-\frac{t}{\Psi(s)}\right)\\
&=\sup_{s\in(0,+\infty)}\left(\frac{C_22^n\Psi^{-1}(t)}{\Psi^{-1}(s)}-\frac{t}{s}\right)\ge C_22^n-1.
\end{align*}
Hence
$$\int_{X\backslash B(x,2^N\Psi^{-1}(t))}p_t(x,y)m(\mathrm{d}y)\le C_1C_{VR}^2\sum_{n=N}^\infty(2^{n+1})^{\alpha_1\vee\alpha_2}\exp(1-C_22^n).$$
Since
$$\sum_{n=0}^\infty(2^{n+1})^{\alpha_1\vee\alpha_2}\exp(1-C_22^n)<+\infty,$$
there exists $N\ge0$ depending only on $C_1$, $C_2$, $C_{VR}$, $\alpha_1$, $\alpha_2$ sufficiently large such that
 $$\int_{X\backslash B(x,2^N\Psi^{-1}(t))}p_t(x,y)m(\mathrm{d}y)\le C_1C_{VR}^2\sum_{n=N}^\infty(2^{n+1})^{\alpha_1\vee\alpha_2}\exp(1-C_22^n)\le\frac{1}{2}.$$
Since $(\mathcal{E},\mathcal{F})$ on $L^2(X;m)$ is conservative, that is, $\int_Xp_t(x,y)m(\mathrm{d}y)=1$, we have
$$\int_{B(x,2^N\Psi^{-1}(t))}p_t(x,y)m(\mathrm{d}y)=1-\int_{X\backslash B(x,2^N\Psi^{-1}(t))}p_t(x,y)m(\mathrm{d}y)\ge\frac{1}{2}.$$
On the other hand
\begin{align*}
&\int_{B(x,2^N\Psi^{-1}(t))}p_t(x,y)m(\mathrm{d}y)\\
&\le\left(\int_{B(x,2^N\Psi^{-1}(t))}p_t(x,y)^2m(\mathrm{d}y)\right)^{1/2}V(x,2^N\Psi^{-1}(t))^{1/2}\\
&\le\left(\int_{X}p_t(x,y)^2m(\mathrm{d}y)\right)^{1/2}V(x,2^N\Psi^{-1}(t))^{1/2}\\
&=p_{2t}(x,x)^{1/2}V(x,2^N\Psi^{-1}(t))^{1/2},
\end{align*}
where we use the Cauchy-Schwarz inequality in the first inequality and the semi-group property in the last equality, hence
$$p_{2t}(x,x)\ge\frac{1}{4V(x,2^N\Psi^{-1}(t))},$$
hence
$$p_{t}(x,x)\ge\frac{1}{4V(x,2^N\Psi^{-1}(\frac{t}{2}))}\ge\frac{\widetilde{{C}}_1}{\Phi(\Psi^{-1}(t))},$$
where $\widetilde{{C}}_1$ is a positive constant depending only on $C_1$, $C_2$, $C_{VR}$, $N$, $\Phi$, $\Psi$.

Second, we prove \ref{eq_NLE}. Indeed, for any $\varepsilon\in(0,\frac{1}{2})$, for any $y\in B(x,\varepsilon\Psi^{-1}(t))$, by \ref{eq_HHK}, we have
\begin{align*}
&|p_t(x,x)-p_t(x,y)|\\
&\le C_{HHK}\frac{\left(\frac{\Psi}{\Phi}\right)(d(x,y))}{t}\\
&\le C_{HHK}\frac{\left(\frac{\Psi}{\Phi}\right)(\varepsilon\Psi^{-1}(t))}{\left(\frac{\Psi}{\Phi}\right)(\Psi^{-1}(t))}\frac{\left(\frac{\Psi}{\Phi}\right)(\Psi^{-1}(t))}{t}\\
&\le C_{HHK}\varepsilon^{\gamma_1\wedge\gamma_2}\frac{1}{\Phi(\Psi^{-1}(t))}.
\end{align*}
There exists $\varepsilon\in(0,\frac{1}{2})$ depending only on $\widetilde{{C}}_1$, $C_{HHK}$, $\gamma_1$, $\gamma_2$ sufficiently small such that
$$C_{HHK}\varepsilon^{\gamma_1\wedge\gamma_2}\le\frac{\widetilde{{C}}_1}{2},$$
hence
$$|p_t(x,x)-p_t(x,y)|\le\frac{\widetilde{{C}}_1}{2\Phi(\Psi^{-1}(t))}\le\frac{1}{2}p_t(x,x).$$
Thus we have
$$p_t(x,y)\ge p_t(x,x)-|p_t(x,x)-p_t(x,y)|\ge p_t(x,x)-\frac{1}{2}p_t(x,x)=\frac{1}{2}p_t(x,x)\ge\frac{\widetilde{{C}}_1}{2\Phi(\Psi^{-1}(t))}$$
for any $y\in B(x,\varepsilon\Psi^{-1}(t))$, that is, we have \ref{eq_NLE}.
\end{proof}

\section{Proof of \texorpdfstring{``\ref{eq_NLE}$\Rightarrow$\ref{eq_HR}"}{"NLE implies HR"}}\label{sec_NLE2HR}

In this section, the strongly recurrent condition (\ref{eq_recurrent}) is intrinsically required to ensure that the Dirichlet form $(\mathcal{E},\mathcal{F})$ on $L^2(X;m)$ indeed induces a resistance form, as will be shown in Lemma \ref{lem_Fe}, after which the arguments rely heavily on Kigami’s resistance form theory \cite{Kig12}.

Recall that a Dirichlet form $(\mathcal{E},\mathcal{F})$ on $L^2(X;m)$ is recurrent if for any non-negative function $f\in L^1(X;m)$, we have either $Gf=+\infty$ or $Gf=0$ $m$-a.e., where $G$ is the Green operator, see \cite[Equation (1.6.2)]{FOT11}. In our setting, the Dirichlet form is recurrent as follows.

\begin{lemma}\label{lem_recurrent}
Let $(X,d,m,\mathcal{E},\mathcal{F})$ be an unbounded MMD space satisfying \ref{eq_VPhi} and \ref{eq_NLE}. Then $(\mathcal{E},\mathcal{F})$ on $L^2(X;m)$ is recurrent.
\end{lemma}

\begin{proof}
Taking any non-zero non-negative function $f\in L^1(X;m)$, we show that
$$Gf(x):=\int_0^{+\infty}\left(\int_Xp_t(x,y)f(y)m(\mathrm{d}y)\right)\mathrm{d}t=+\infty\text{ for }m\text{-a.e. }x\in X.$$
Indeed, by \ref{eq_VPhi} and \ref{eq_NLE}, we have
\begin{align*}
&\int_0^{+\infty}p_t(x,y)\mathrm{d}t\ge\int_{\Psi(\frac{d(x,y)}{\varepsilon})}^{+\infty}p_t(x,y)\mathrm{d}t\ge\int_{\Psi(\frac{d(x,y)}{\varepsilon})}^{+\infty}\frac{C}{V(x,\Psi^{-1}(t))}\mathrm{d}t\\
&\asymp\int_{\Psi(\frac{d(x,y)}{\varepsilon})}^{+\infty}\frac{1}{\Phi(\Psi^{-1}(t))}\mathrm{d}t\ge\int_{\Psi(\frac{d(x,y)}{\varepsilon})\vee1}^{+\infty}\frac{\mathrm{d}t}{t^{\alpha_2/\beta_2}}=+\infty.
\end{align*}
Since $f$ is non-zero non-negative, there exists $a>0$ such that $m(\{f>a\})>0$, then
$$Gf(x)=\int_X\left(\int_0^{+\infty}p_t(x,y)\mathrm{d}t\right)f(y)m(\mathrm{d}y)\ge a\int_{\{f>a\}}\left(\int_0^{+\infty}p_t(x,y)\mathrm{d}t\right)m(\mathrm{d}y)=+\infty$$
for $m$-a.e. $x\in X$. Hence $(\mathcal{E},\mathcal{F})$ on $L^2(X;m)$ is recurrent.
\end{proof}

To begin with, let us introduce effective resistances. Let $A,B$ be two non-empty subsets of $X$. The effective resistance between $A$ and $B$ is defined as follows.
$$R(A,B)=\left(\inf\left\{\mathcal{E}(u,u):u\in\mathcal{F},u=0\text{ on }A,u=1\text{ on }B\right\}\right)^{-1},$$
here we use the convention that $\inf\emptyset=+\infty$, $0^{-1}=+\infty$  and $(+\infty)^{-1}=0$. For $x,y\in X$, write
$$R(x,A)=R(\{x\},A),R(x,y)=R(\{x\},\{y\}).$$
By definition, it is obvious that
$$R(x,y)=\sup\left\{\frac{|u(x)-u(y)|^2}{\mathcal{E}(u,u)}:u\in\mathcal{F},\mathcal{E}(u,u)>0\right\},$$
which implies that
$$|u(x)-u(y)|^2\le R(x,y)\mathcal{E}(u,u)\text{ for any }u\in\mathcal{F},x,y\in X,$$
and if $A_1\subseteq A_2$, $B_1\subseteq B_2$, then
$$R(A_1,B_1)\ge R(A_2,B_2).$$

\begin{proposition}\label{prop_res}
Let $(X,d,m,\mathcal{E},\mathcal{F})$ be an unbounded MMD space satisfying \ref{eq_VPhi},\\
\noindent \ref{eq_UHK} and \ref{eq_NLE}. Then we have two-sided resistance estimates \ref{eq_res} as follows. There exists $C_R\in(0,+\infty)$ such that for any $x,y\in X$ with $x\ne y$, we have
\begin{align*}\label{eq_res}\tag*{R($\Phi,\Psi$)}
\frac{1}{C_R}\frac{\Psi(d(x,y))}{\Phi(d(x,y))}\le R(x,y)\le C_R\frac{\Psi(d(x,y))}{\Phi(d(x,y))}.
\end{align*}
\end{proposition}

\begin{remark}
In the strongly recurrent setting, two-sided heat kernel estimates and two-sided resistance estimates are indeed equivalent. This equivalence was proved in \cite{BCK05} for graphs, in \cite{Kum04} for graphs and resistance forms, and in \cite{Hu08} for metric spaces using an analytical approach.
\end{remark}

Under the assumptions of the above result, by Proposition \ref{prop_HK}, we have \ref{eq_PI} and \ref{eq_CS}.

The following result will play an important role in this section.

\begin{lemma}[Morrey-Sobolev inequality]\label{lem_MS}
Let $(X,d,m,\mathcal{E},\mathcal{F})$ be an unbounded MMD space satisfying \ref{eq_VPhi} and \ref{eq_PI}. Then we have the following Morrey-Sobolev inequality \ref{eq_MS}. There exists $C_{MS}\in(0,+\infty)$ such that for any $u\in\mathcal{F}$, for $m$-a.e. $x,y\in X$ with $x\ne y$, we have
\begin{align*}\label{eq_MS}\tag*{MS($\Phi,\Psi$)}
|u(x)-u(y)|^2\le C_{MS}\frac{\Psi(d(x,y))}{\Phi(d(x,y))}\mathcal{E}(u,u).
\end{align*}
Hence any function in $\mathcal{F}$ has a continuous version, or equivalently, $\mathcal{F}\subseteq C(X)$.
\end{lemma}

The proof is standard using a telescopic technique, see also \cite[Page 1654]{BCK05}.

\begin{proof}
Let $x,y$ be two different Lebesgue points of $u\in\mathcal{F}\subseteq L^2(X;m)$. Denote $r=d(x,y)$. Then
$$|u(x)-u(y)|\le|u(x)-u_{B(x,r)}|+|u_{B(x,r)}-u_{B(y,r)}|+|u(y)-u_{B(y,r)}|.$$

For any $n\in \mathbb{Z}$, denote $B_n=B(x,2^{-n}r)$. By \ref{eq_VPhi} and \ref{eq_PI}, we have
\begin{align*}
&|u_{B_{n+1}}-u_{B_n}|\le\dashint_{B_{n+1}}|u-u_{B_n}|\mathrm{d}m\lesssim\dashint_{B_n}|u-u_{B_n}|\mathrm{d}m\le\left(\dashint_{B_n}|u-u_{B_n}|^2\mathrm{d}m\right)^{1/2}\\
&\lesssim\left(\frac{1}{\Phi(2^{-n}r)}\Psi(2^{-n}r)\int_{B_{n-1}}\mathrm{d}\Gamma(u,u)\right)^{1/2}\le\left(\frac{\Psi(2^{-n}r)}{\Phi(2^{-n}r)}\mathcal{E}(u,u)\right)^{1/2},
\end{align*}
hence
\begin{align*}
&|u(x)-u_{B(x,r)}|=|u(x)-u_{B_0}|=\lim_{n\to+\infty}|u_{B_{n}}-u_{B_0}|\le\varliminf_{n\to+\infty}\sum_{k=0}^{n-1}|u_{B_{k+1}}-u_{B_k}|\\
&=\sum_{n=0}^\infty|u_{B_{n+1}}-u_{B_n}|\lesssim\sum_{n=0}^\infty\left(\frac{\Psi(2^{-n}r)}{\Phi(2^{-n}r)}\mathcal{E}(u,u)\right)^{1/2}\\
&=\left(\frac{\Psi(r)}{\Phi(r)}\mathcal{E}(u,u)\right)^{1/2}\sum_{n=0}^\infty\left(\frac{\left(\frac{\Psi}{\Phi}\right)(2^{-n}r)}{\left(\frac{\Psi}{\Phi}\right)(r)}\right)^{1/2}\le\left(\frac{\Psi(r)}{\Phi(r)}\mathcal{E}(u,u)\right)^{1/2}\sum_{n=0}^{\infty}2^{-\frac{(\gamma_1\wedge\gamma_2)n}{2}}\\
&\asymp\left(\frac{\Psi(r)}{\Phi(r)}\mathcal{E}(u,u)\right)^{1/2}.
\end{align*}
Similarly, we have
$$|u(y)-u_{B(y,r)}|\lesssim\left(\frac{\Psi(r)}{\Phi(r)}\mathcal{E}(u,u)\right)^{1/2}.$$
Note that
\begin{align*}
&|u_{B(x,r)}-u_{B(y,r)}|\le\dashint_{B(x,r)}\dashint_{B(y,r)}|u(z_1)-u(z_2)|m(\mathrm{d}z_1)m(\mathrm{d}z_2)\\
&\le\left(\dashint_{B(x,r)}\dashint_{B(y,r)}|u(z_1)-u(z_2)|^2m(\mathrm{d}z_1)m(\mathrm{d}z_2)\right)^{1/2}\\
&\asymp\frac{1}{\Phi(r)}\left(\int_{B(x,r)}\int_{B(y,r)}|u(z_1)-u(z_2)|^2m(\mathrm{d}z_1)m(\mathrm{d}z_2)\right)^{1/2}\\
&\le\frac{1}{\Phi(r)}\left(\int_{B(x,2r)}\int_{B(x,2r)}|u(z_1)-u(z_2)|^2m(\mathrm{d}z_1)m(\mathrm{d}z_2)\right)^{1/2}\\
&=\frac{1}{\Phi(r)}\left(2V(x,2r)\int_{B(x,2r)}\left(u-u_{B(x,2r)}\right)^2\mathrm{d}m\right)^{1/2}\\
&\lesssim\frac{1}{\Phi(r)^{1/2}}\left(\Psi(2r)\int_{B(x,4r)}\mathrm{d}\Gamma(u,u)\right)\lesssim\left(\frac{\Psi(r)}{\Phi(r)}\mathcal{E}(u,u)\right)^{1/2}.
\end{align*}

In summary, we have
$$|u(x)-u(y)|\lesssim\left(\frac{\Psi(r)}{\Phi(r)}\mathcal{E}(u,u)\right)^{1/2}.$$
\end{proof}

\begin{proof}[Proof of Proposition \ref{prop_res}]
The upper bound is obvious by \ref{eq_MS}. We use the idea of an analytical proof of \cite[Theorem 3.4]{Hu08} to give the lower bound. By the spectral calculus, it is easy to see that
$$\mathcal{E}(P_tu,P_tu)\le\frac{1}{2et}(u,u)\text{ for any }t\in(0,+\infty),u\in L^2(X;m),$$
where $\{P_t\}$ is the heat semi-group corresponding to the Dirichlet form $(\mathcal{E},\mathcal{F})$ on $L^2(X;m)$.

For any $t\in(0,+\infty)$, $x\in X$, take $u=p_t(x,\cdot)\in L^2(X;m)$. By the semi-group property, we have $P_tu=p_{2t}(x,\cdot)$ and $(u,u)=p_{2t}(x,x)$, hence
$$\mathcal{E}(p_{2t}(x,\cdot),p_{2t}(x,\cdot))\le\frac{1}{2et}p_{2t}(x,x),$$
clearly
$$\mathcal{E}(p_t(x,\cdot),p_t(x,\cdot))\le\frac{1}{et}p_t(x,x).$$
For any $y\in X$, we have
\begin{equation}\label{eq_res1}
p_t(x,x)-p_t(x,y)\le R(x,y)^{1/2}\mathcal{E}(p_t(x,\cdot),p_t(x,\cdot))^{1/2}\le R(x,y)^{1/2}\left(\frac{1}{et}p_t(x,x)\right)^{1/2}.
\end{equation}
By \ref{eq_NLE}, we have
$$p_t(x,x)\ge\frac{C_1}{V(x,\Psi^{-1}(t))}.$$
By \ref{eq_UHK}, we have
$$p_t(x,y)\le\frac{C_2}{V(x,\Psi^{-1}(t))}\exp\left(-\Upsilon(C_3d(x,y),t)\right),$$
and
$$p_t(x,x)\le\frac{C_2}{V(x,\Psi^{-1}(t))}.$$
Hence
\begin{align*}
&p_t(x,x)-p_t(x,y)\\
&\ge\frac{C_1}{V(x,\Psi^{-1}(t))}-\frac{C_2}{V(x,\Psi^{-1}(t))}\exp\left(-\Upsilon(C_3d(x,y),t)\right)\\
&=\frac{1}{V(x,\Psi^{-1}(t))}\left(C_1-C_2\exp\left(-\Upsilon(C_3d(x,y),t)\right)\right),
\end{align*}
where
\begin{align*}
&\Upsilon(C_3d(x,y),t)=\sup_{s\in(0,+\infty)}\left(\frac{C_3d(x,y)}{s}-\frac{t}{\Psi(s)}\right)\\
&=\sup_{s\in(0,+\infty)}\left(\frac{C_3d(x,y)}{\Psi^{-1}(s)}-\frac{t}{s}\right)\ge \frac{C_3d(x,y)}{\Psi^{-1}(t)}-1.
\end{align*}
Take $\varepsilon\in(0,1)$, let $t=\Psi(\varepsilon d(x,y))$, then
\begin{equation}\label{eq_res2}
p_t(x,x)-p_t(x,y)\ge\frac{1}{V(x,\Psi^{-1}(t))}\left(C_1-C_2\exp\left(1-\frac{C_3}{\varepsilon}\right)\right)\ge\frac{C_1}{2V(x,\Psi^{-1}(t))},
\end{equation}
where $\varepsilon\in(0,1)$ depending only on $C_1$, $C_2$, $C_3$ is sufficiently small such that
$$C_2\exp\left(1-\frac{C_3}{\varepsilon}\right)\le\frac{C_1}{2}.$$
Plugging Equation (\ref{eq_res2}) into Equation (\ref{eq_res1}), we have
$$\frac{C_1}{2V(x,\Psi^{-1}(t))}\le R(x,y)^{1/2}\left(\frac{1}{et}p_t(x,x)\right)^{1/2}\le R(x,y)^{1/2}\left(\frac{1}{et}\frac{C_2}{V(x,\Psi^{-1}(t))}\right)^{1/2},$$
therefore
\begin{align*}
&R(x,y)^{1/2}\ge \frac{C_1}{2V(x,\Psi^{-1}(t))}\left(\frac{etV(x,\Psi^{-1}(t))}{C_2}\right)^{1/2}\\
&\asymp\left(\frac{t}{\Phi(\Psi^{-1}(t))}\right)^{1/2}=\left(\frac{\Psi(\varepsilon d(x,y))}{\Phi(\varepsilon d(x,y))}\right)^{1/2}\asymp\left(\frac{\Psi(d(x,y))}{\Phi(d(x,y))}\right)^{1/2},
\end{align*}
that is,
$$R(x,y)\gtrsim \frac{\Psi(d(x,y))}{\Phi(d(x,y))}.$$
\end{proof}

Let us recall the definition of resistance forms as follows.

\begin{definition}[{\cite[DEFINITION 3.1]{Kig12}}]\label{def_RF}
Let $X$ be a set. A pair $(\mathcal{E},\mathcal{G})$ is called a resistance form on $X$ if the following five conditions are satisfied.
\begin{enumerate}[label=(RF\arabic*),ref=(RF\arabic*)]
\item\label{def_RF1} $\mathcal{G}$ is a linear subspace of $l(X)$ containing constant functions, where $l(X)$ is the space of all real-valued functions on $X$. $\mathcal{E}$ is a non-negative symmetric quadratic form on $\mathcal{G}$. For $u\in\mathcal{G}$, $\mathcal{E}(u,u)=0$ if and only if $u$ is a constant function.
\item\label{def_RF2} Let $\sim$ be the equivalence relation on $\mathcal{G}$ given by $u,v\in\mathcal{G}$, $u\sim v$ if and only if $u-v$ is a constant function. Then $(\mathcal{G}/\!\sim,\mathcal{E})$ is a Hilbert space.
\item\label{def_RF3} For any $x,y\in X$ satisfying $x\ne y$, there exists $u\in\mathcal{G}$ such that $u(x)\ne u(y)$.
\item\label{def_RF4} For any $x,y\in X$, we have
$$\sup\left\{\frac{|u(x)-u(y)|^2}{\mathcal{E}(u,u)}:u\in\mathcal{G},\mathcal{E}(u,u)>0\right\}$$
is finite, where we use the convention that $\sup\emptyset=0$.
\item\label{def_RF5} For any $u\in\mathcal{G}$, let $\overline{{u}}=(u\vee0)\wedge1$, then $\overline{{u}}\in\mathcal{G}$ and $\mathcal{E}(\overline{{u}},\overline{{u}})\le\mathcal{E}(u,u)$.
\end{enumerate}
\end{definition}

Recall that for a Dirichlet form $(\mathcal{E},\mathcal{F})$ on $L^2(X;m)$, its extended Dirichlet space is defined as the family of all $m$-measurable functions $u$ which is finite $m$-a.e. and there exists an $\mathcal{E}$-Cauchy sequence $\left\{u_n\right\}\subseteq\mathcal{F}$ such that $u_n\to u$ $m$-a.e.. The sequence $\left\{u_n\right\}\subseteq\mathcal{F}$ is called an approximation sequence of $u\in\mathcal{F}_e$. We can extend $\mathcal{E}$ to $\mathcal{F}_e$ and $\mathcal{F}=\mathcal{F}_e\cap L^2(X;m)$. See \cite[Section 1.5]{FOT11} for more details. We can also define the harmonicity of a function in the extended Dirichlet space as follows. Let $D$ be a non-empty open subset of $X$. We say that $u\in\mathcal{F}_e$ is harmonic in $D$ if $\mathcal{E}(u,\varphi)=0$ for any $\varphi\in\mathcal{F}\cap C_c(D)$.

We collect some results about extended Dirichlet spaces as follows.

\begin{lemma}\label{lem_Fe}
Let $(X,d,m,\mathcal{E},\mathcal{F})$ be an unbounded MMD space satisfying \ref{eq_VPhi}, \ref{eq_UHK} and \ref{eq_NLE}.
\begin{enumerate}[label=(\arabic*),ref=(\arabic*)]
\item\label{lem_Fe_MS} The Morrey-Sobolev inequality \ref{eq_MS} holds for any $u\in\mathcal{F}_e$, hence $\mathcal{F}_e\subseteq C(X)$.
\item\label{lem_Fe_RF} $(\mathcal{E},\mathcal{F}_e)$ is a resistance form on $X$.
\item\label{lem_Fe_res} For any $x,y\in X$, we have
$$R(x,y)=\sup\left\{\frac{|u(x)-u(y)|^2}{\mathcal{E}(u,u)}:u\in\mathcal{F}_e,\mathcal{E}(u,u)>0\right\}.$$
\item\label{lem_Fe_har} Let $u\in\mathcal{F}_e$ be harmonic in a non-empty open subset $D$ of $X$. Then $\mathcal{E}(u,\varphi)=0$ for any $\varphi\in\mathcal{F}_{e,D}$, where
\begin{equation}\label{eq_FeD}
\mathcal{F}_{e,D}=\left\{u\in\mathcal{F}_e\subseteq C(X):u=0\text{ on }X\backslash D\right\}.
\end{equation}
\end{enumerate}
\end{lemma}

\begin{proof}
\ref{lem_Fe_MS} It is obvious by the definition of $\mathcal{F}_e$ and Lemma \ref{lem_MS}.

\ref{lem_Fe_RF} Since the Dirichlet form $(\mathcal{E},\mathcal{F})$ on $L^2(X;m)$ is recurrent, by \cite[Theorem 1.6.3]{FOT11}, we have $1\in\mathcal{F}_e$ and $\mathcal{E}(1,1)=0$. For $u\in\mathcal{F}_e$, if $\mathcal{E}(u,u)=0$, then by \ref{eq_MS}, we have $u$ is a constant function. Hence we have \ref{def_RF1}.

By \ref{def_RF1}, we have $(\mathcal{F}_e/\!\sim,\mathcal{E})$ is an inner product space. We only need to prove the completeness. Indeed, let $\{u_n\}\subseteq\mathcal{F}_e/\!\sim$ be an $\mathcal{E}$-Cauchy sequence, then $\{u_n\}$ is $\mathcal{E}$-bounded. Fix arbitrary $x_0\in X$. By replacing $u_n$ by $u_n-u_n(x_0)$ in $\mathcal{F}_e/\!\sim$, we may assume that $u_n(x_0)=0$ for any $n\ge1$.

For any $N\ge1$, for any $x\in B(x_0,N)$, for any $n\ge1$, by \ref{eq_MS}, we have
$$|u_n(x)|^2=|u_n(x)-u_n(x_0)|^2\le C_{MS}\frac{\Psi(d(x,x_0))}{\Phi(d(x,x_0))}\mathcal{E}(u_n,u_n)\le C_{MS}\frac{\Psi(N)}{\Phi(N)}\sup_{n\ge1}\mathcal{E}(u_n,u_n).$$
By \ref{eq_MS} again, using the Arzel\`a-Ascoli lemma and the diagonal argument, there exists a subsequence, still denoted by $\{u_n\}$, and $u\in C(X)$ such that $u_n$ converges uniformly to $u$ on any relatively compact open subset of $X$. By the definition of $\mathcal{F}_e$, we have $u\in\mathcal{F}_e$ and $\{u_n\}$ is $\mathcal{E}$-convergent to $u$. Hence $(\mathcal{F}_e/\!\sim,\mathcal{E})$ is complete, then we have \ref{def_RF2}.

\ref{def_RF3} follows from the regularity of the Dirichlet form $(\mathcal{E},\mathcal{F})$ on $L^2(X;m)$.

\ref{def_RF4} follows from \ref{eq_MS}.

\ref{def_RF5} follows from \cite[Corollary 1.6.3]{FOT11}.

\ref{lem_Fe_res} We need to prove that
$$\sup\left\{\frac{|u(x)-u(y)|^2}{\mathcal{E}(u,u)}:u\in\mathcal{F},\mathcal{E}(u,u)>0\right\}=\sup\left\{\frac{|u(x)-u(y)|^2}{\mathcal{E}(u,u)}:u\in\mathcal{F}_e,\mathcal{E}(u,u)>0\right\}.$$
Indeed, ``$\le$" is trivial by $\mathcal{F}\subseteq\mathcal{F}_e$. ``$\ge$": For any $u\in\mathcal{F}_e$ with $\mathcal{E}(u,u)>0$, there exists $\{u_n\}\subseteq\mathcal{F}$ which is $\mathcal{E}$-Cauchy and converges to $u$ $m$-a.e.. Since $\{u_n\}$ is $\mathcal{E}$-convergent to $u$, we may assume that $\mathcal{E}(u_n,u_n)>0$ for any $n\ge1$. Take $x_0\in X$ such that $\lim_{n\to+\infty}u_n(x_0)=u(x_0)\in \mathbb{R}$, then $\{u_n(x_0)\}$ is bounded. For any $N\ge1$, for any $x\in B(x_0,N)$, for any $n\ge1$, by \ref{eq_MS}, we have
\begin{align*}
&|u_n(x)|^2\le2\left(|u_n(x)-u_n(x_0)|^2+|u_n(x_0)|^2\right)\\
&\le2\left(C_{MS}\frac{\Psi(d(x,x_0))}{\Phi(d(x,x_0))}\mathcal{E}(u_n,u_n)+\sup_{n\ge1}|u_n(x_0)|^2\right)\\
&\le2\left(C_{MS}\frac{\Psi(N)}{\Phi(N)}\sup_{n\ge1}\mathcal{E}(u_n,u_n)+\sup_{n\ge1}|u_n(x_0)|^2\right).
\end{align*}
By \ref{eq_MS} again, using the Arzel\`a-Ascoli lemma and the diagonal argument, there exists a subsequence, still denoted by $\{u_n\}$, and $v\in C(X)$ such that $u_n$ converges uniformly to $v$ on any relatively compact open subset of $X$. Since $u_n$ also converges to $u$ $m$-a.e. and $u,v\in C(X)$, we have $u=v$. Hence
$$\frac{|u(x)-u(y)|^2}{\mathcal{E}(u,u)}=\lim_{n\to+\infty}\frac{|u_n(x)-u_n(y)|^2}{\mathcal{E}(u_n,u_n)}\le\text{LHS}.$$
Taking supremum with respect to $u\in\mathcal{F}_e$ with $\mathcal{E}(u,u)>0$, we have ``$\ge$".

\ref{lem_Fe_har} Since $(\mathcal{E},\mathcal{F}_e)$ is a resistance form on $X$, by \cite[PROPOSITION 9.13]{Kig12}, we have the continuity and the quasi-continuity coincide. By \cite[Equation (2.3.14)]{FOT11}, we have
$$\mathcal{F}_{e,D}=\left\{u\in\mathcal{F}_e:\widetilde{{u}}=0\text{ q.e. on }X\backslash D\right\},$$
where $\widetilde{{u}}$ is a quasi-continuous modification of $u\in\mathcal{F}_e$.

Since the Dirichlet form $(\mathcal{E},\mathcal{F})$ on $L^2(X;m)$ is regular, by \cite[Exercise 1.4.1]{FOT11}, we have $\mathcal{F}\cap C_c(X)$ is a special standard core of $\mathcal{E}$. Then \ref{lem_Fe_har} follows from the equivalence in \cite[Theorem 2.3.3 (ii)]{FOT11}.
\end{proof}

By resistance form theory, we can solve ``boundary value problems" and take traces of resistance forms as follows.

\begin{lemma}[{\cite[Section 8]{Kig12}}]\label{lem_trace}
Let $(X,d,m,\mathcal{E},\mathcal{F})$ be an unbounded MMD space satisfying \ref{eq_VPhi}, \ref{eq_UHK} and \ref{eq_NLE}. Let $Y\subsetneqq X$ be a non-empty closed subset and $\mathcal{F}_e|_Y=\{u|_{Y}:u\in\mathcal{F}_e\}$. Then for any $u\in\mathcal{F}_e$, there exists a unique $v\in\mathcal{F}_e$ which is harmonic in $X\backslash Y$ such that $v=u$ on $Y$. Moreover, $v$ is the unique function in $\mathcal{F}_e$ that minimizes the following variational problem
$$\inf\left\{\mathcal{E}(v,v):v\in\mathcal{F}_e,v=u\text{ on }Y\right\}.$$
Let $h_Y:\mathcal{F}_e|_Y\to\mathcal{F}_e$ be given by $u|_{Y}\mapsto v$ for $u\in\mathcal{F}_e$ as above. Let
$$\mathcal{E}|_Y(u,v)=\mathcal{E}(h_Y(u),h_Y(v))$$
for $u,v\in\mathcal{F}_e|_Y$. Then $(\mathcal{E}|_Y,\mathcal{F}_e|_Y)$ is a resistance form on $Y$.
\end{lemma}

We have the following maximum principle.

\begin{corollary}[Maximum principle]\label{cor_max}
Let $(X,d,m,\mathcal{E},\mathcal{F})$ be an unbounded MMD space satisfying \ref{eq_VPhi}, \ref{eq_UHK} and \ref{eq_NLE}. Let $D\subsetneqq X$ be a non-empty open subset of $X$. Let $u,v\in\mathcal{F}_e$ be harmonic in $D$ and satisfy $u\le v$ on $X\backslash D$, then $u\le v$ in $D$. In particular, $\inf_{X\backslash D}u\le u\le\sup_{X\backslash D}u$.
\end{corollary}

\begin{proof}
We only need to show that if $u\in\mathcal{F}_e$ is harmonic in $D$ and satisfies $u\le0$ on $X\backslash D$, then $u\le 0$ in $D$. Indeed, $u\wedge0\in\mathcal{F}_e$ satisfies $u\wedge0=u$ on $X\backslash D$, by Lemma \ref{lem_trace}, we have
$$\mathcal{E}(u,u)\ge\mathcal{E}(u\wedge0,u\wedge0)\ge\inf\left\{\mathcal{E}(v,v):v\in\mathcal{F}_e,v=u\text{ on }X\backslash D\right\}=\mathcal{E}(u,u).$$
Hence the above inequalities are indeed equalities. By the uniqueness of $u\in\mathcal{F}_e$ from Lemma \ref{lem_trace}, we have
$u\wedge0=u$, hence $u\le0$ in $D$.
\end{proof}

For certain pairs of sets, the effective resistance can be realized by the energy of a suitable function, as follows.

\begin{corollary}\label{cor_attain}
Let $(X,d,m,\mathcal{E},\mathcal{F})$ be an unbounded MMD space satisfying \ref{eq_VPhi}, \ref{eq_UHK} and \ref{eq_NLE}. Let $A,B$ be two non-empty closed subsets of $X$ satisfying $A\cap B=\emptyset$ and $A\cup B\subsetneqq X$. Assume that one of the following conditions holds:
\begin{enumerate}[label=(\alph*),ref=(\alph*)]
\item\label{item_attain_a} $X\backslash A$ is relatively compact and $B$ is compact.
\item\label{item_attain_b} $A,B$ are compact.
\end{enumerate}
Then there exists a unique $v\in\mathcal{F}_e$ with $v=0$ on $A$ and $v=1$ on $B$ such that
$$R(A,B)^{-1}=\mathcal{E}(v,v)=\inf\left\{\mathcal{E}(v,v):v\in\mathcal{F}_e,v=0\text{ on }A,v=1\text{ on }B\right\}.$$
Moreover, $1_A,1_B\in\mathcal{F}_e|_{A\cup B}$, $v=h_{A\cup B}(1_B)\in\mathcal{F}_e$ and
$$R(A,B)^{-1}=\mathcal{E}(h_{A\cup B}(1_A),h_{A\cup B}(1_A))=\mathcal{E}(h_{A\cup B}(1_B),h_{A\cup B}(1_B)).$$
\end{corollary}

\begin{proof}
Since $B$ is compact, $X\backslash A$ is open and $B\subseteq X\backslash A$, by the regularity of $(\mathcal{E},\mathcal{F})$ on $L^2(X;m)$, there exists $u\in\mathcal{F}\cap C_c(X)\subseteq\mathcal{F}_e$ such that $0\le u\le 1$ in $X$, $u=1$ on $B$ and $\mathrm{supp}(u)\subseteq X\backslash A$. Hence $u|_{A\cup B}=1_B\in\mathcal{F}_e|_{A\cup B}$, $(1-u)|_{A\cup B}=1_A\in\mathcal{F}_e|_{A\cup B}$.

Since $A\cup B\subsetneqq X$ is a non-empty closed subset, by Lemma \ref{lem_trace}, there exists a unique $v=h_{A\cup B}(u|_{A\cup B})=h_{A\cup B}(1_B)\in\mathcal{F}_e$ which is harmonic in $X\backslash(A\cup B)$ and $v=u$ on $A\cup B$, that is, $v|_{A\cup B}=1_B$, or $v=0$ on $A$ and $v=1$ on $B$. Moreover, $v\in\mathcal{F}_e$ also minimizes the variational problem
$$\inf\left\{\mathcal{E}(v,v):v\in\mathcal{F}_e,v=0\text{ on }A,v=1\text{ on }B\right\},$$
and
\begin{align*}
&\mathcal{E}(v,v)=\mathcal{E}|_{A\cup B}(1_B,1_B)=\mathcal{E}|_{A\cup B}(1_A,1_A)\\
&=\mathcal{E}(h_{A\cup B}(1_B),h_{A\cup B}(1_B))=\mathcal{E}(h_{A\cup B}(1_A),h_{A\cup B}(1_A)).
\end{align*}

It remains to prove that if either (a) or (b) holds, then 
$$\mathcal{E}(v,v)=\inf\left\{\mathcal{E}(v,v):v\in\mathcal{F},v=0\text{ on }A,v=1\text{ on }B\right\},$$
where the RHS is equal to $R(A,B)^{-1}$ by definition. Since $\mathcal{F}\subseteq\mathcal{F}_e$, we only need to show ``$\ge$". By the maximum principle (Corollary \ref{cor_max}), we have $0\le v\le1$ in $X$.

Assuming \ref{item_attain_a}, since $0\le v\le 1$ in $X$, $v=0$ on $A$ and $X\backslash A$ is a relatively compact open set, we have $v\in\mathcal{F}_e\cap L^2(X;m)=\mathcal{F}$, ``$\ge$" is obvious.

Assuming \ref{item_attain_b}, since $A,B$ are both compact, by a similar proof to Lemma \ref{lem_Fe} \ref{lem_Fe_RF}, there exists an approximation  sequence $\{v_n\}\subseteq\mathcal{F}\subseteq C(X)$ of $v\in\mathcal{F}_e$ that converges uniformly to $v$ on $A\cup B$. For any $\varepsilon\in(0,1/4)$, there exists $N\ge1$ such that for any $n\ge N$, we have $|v_n-v|<\varepsilon$ on $A\cup B$, hence $-\varepsilon<v_n<\varepsilon$ on $A$ and $1-\varepsilon<v_n<1+\varepsilon$ on $B$. Let $\varphi:\mathbb{R}\to \mathbb{R}$ be given by
$$\varphi(t)=\frac{\frac{1}{2}}{\frac{1}{2}-\varepsilon}\left((t\wedge(1-\varepsilon))\vee\varepsilon-\frac{1}{2}\right)+\frac{1}{2},t\in \mathbb{R},$$
that is, $\varphi\in C(\mathbb{R})$ satisfies $\varphi=0$ on $(-\infty,\varepsilon]$, $\varphi=1$ on $[1-\varepsilon,+\infty)$ and $\varphi$ is linear on $[\varepsilon,1-\varepsilon]$. Then $\varphi(0)=0$ and $|\varphi(t)-\varphi(s)|\le\frac{\frac{1}{2}}{\frac{1}{2}-\varepsilon}|t-s|$ for any $t,s\in \mathbb{R}$. For any $n\ge N$, let $\widetilde{{v}}_n=\varphi(v_n)$, then $\widetilde{{v}}_n\in\mathcal{F}$ satisfies $\widetilde{{v}}_n=0$ on $A$, $\widetilde{{v}}_n=1$ on $B$ and $\mathcal{E}(\widetilde{{v}}_n,\widetilde{{v}}_n)\le\left(\frac{\frac{1}{2}}{\frac{1}{2}-\varepsilon}\right)^2\mathcal{E}(v_n,v_n)$, hence
$$\inf\left\{\mathcal{E}(v,v):v\in\mathcal{F},v=0\text{ on }A,v=1\text{ on }B\right\}\le\mathcal{E}(\widetilde{{v}}_n,\widetilde{{v}}_n)\le\left(\frac{\frac{1}{2}}{\frac{1}{2}-\varepsilon}\right)^2\mathcal{E}(v_n,v_n).$$
Letting $n\to+\infty$, we have
$$\inf\left\{\mathcal{E}(v,v):v\in\mathcal{F},v=0\text{ on }A,v=1\text{ on }B\right\}\le\left(\frac{\frac{1}{2}}{\frac{1}{2}-\varepsilon}\right)^2\mathcal{E}(v,v).$$
Since $\varepsilon\in(0,\frac{1}{4})$ is arbitrary, we have
$$\inf\left\{\mathcal{E}(v,v):v\in\mathcal{F},v=0\text{ on }A,v=1\text{ on }B\right\}\le\mathcal{E}(v,v),$$
that is, we have ``$\ge$".
\end{proof}

\begin{remark}
We will be interested in the case $A=X\backslash B(x,r)$, $B=\{y\}$ in \ref{item_attain_a}, and in the case $A=\{x\}$, $B=\{y\}$ in \ref{item_attain_b}, where $x,y\in X$, $r\in(0,+\infty)$.
\end{remark}

\begin{lemma}\label{lem_ptball}
Let $(X,d,m,\mathcal{E},\mathcal{F})$ be an unbounded MMD space satisfying \ref{eq_VPhi}, \ref{eq_UHK} and \ref{eq_NLE}. Then there exists $C\in(0,+\infty)$ such that for any $x_0\in X$, for any $r\in(0,+\infty)$, we have
$$R(x_0,X\backslash B(x_0,r))\ge C\frac{\Psi(r)}{\Phi(r)}.$$
\end{lemma}

The proof is essentially the same as the proof of \cite[LEMMA 2.4]{BCK05}, \cite[Lemma 4.1]{Kum04} and \cite[Proposition 5.3]{Hu08}.

\begin{proof}
Denote $B=B(x_0,r)$. For any $x\in B\backslash\left(\frac{1}{2}B\right)$, by Corollary \ref{cor_attain}, there exists $h_x\in\mathcal{F}_e$ with $h_x(x_0)=1$, $h_x(x)=0$ and $0\le h_x\le 1$ in $X$ such that
$$\mathcal{E}(h_x,h_x)=R(x_0,x)^{-1}\le\left(\frac{1}{C_R}\frac{\Psi(d(x_0,x))}{\Phi(d(x_0,x))}\right)^{-1}=C_R\frac{\Phi(d(x_0,x))}{\Psi(d(x_0,x))}\le C_R\frac{\Phi(\frac{r}{2})}{\Psi(\frac{r}{2})},$$
where $C_R$ is the constant in \ref{eq_res}. Take $\varepsilon\in(0,\frac{1}{4})$ to be determined. For any $y\in B(x,\varepsilon r)$, we have
\begin{align*}
&|h_x(y)|=|h_x(x)-h_x(y)|\le R(x,y)^{1/2}\mathcal{E}(h_x,h_x)^{1/2}\\
&\le\left(C_R\frac{\Psi(d(x,y))}{\Phi(d(x,y))}C_R\frac{\Phi(\frac{r}{2})}{\Psi(\frac{r}{2})}\right)^{1/2}\le C_R\left(\frac{\frac{\Psi(\varepsilon r)}{\Phi(\varepsilon r)}}{\frac{\Psi(\frac{r}{2})}{\Phi(\frac{r}{2})}}\right)^{1/2}\le C_R(2\varepsilon)^{\frac{{\gamma_1\wedge\gamma_2}}{2}}.
\end{align*}
Then there exists $\varepsilon\in(0,\frac{1}{4})$ depending only on $C_R,\gamma_1,\gamma_2$ sufficiently small such that
$$|h_x(y)|\le C_R(2\varepsilon)^{\frac{{\gamma_1\wedge\gamma_2}}{2}}\le\frac{1}{2}\text{ for any }y\in B(x,\varepsilon r).$$
By \ref{eq_VPhi}, there exists $N\ge1$ depending only on $C_{VR},\alpha_1,\alpha_2, \varepsilon$, hence only on $C_{VR}$, $C_R$, $\alpha_1$, $\alpha_2$, $\beta_1$, $\beta_2$, such that $B\backslash\left(\frac{1}{2}B\right)$ can be covered by at most $N$ balls  with center in $B\backslash\left(\frac{1}{2}B\right)$ and radius $\varepsilon r$. Take $x_1,\ldots,x_N\in B\backslash\left(\frac{1}{2}B\right)$ such that $B\backslash\left(\frac{1}{2}B\right)\subseteq\cup_{1\le i\le N}B(x_i,\varepsilon r)$, let $h=\min_{1\le i\le N}h_{x_i}$. Then $h\in\mathcal{F}_e$ satisfies $h(x_0)=1$, $0\le h\le 1$ in $X$, $h\le\frac{1}{2}$ in $B\backslash\left(\frac{1}{2}B\right)$ and
$$\mathcal{E}(h,h)\le\sum_{i=1}^N\mathcal{E}(h_{x_i},h_{x_i})\le \sum_{i=1}^NC_R\frac{\Phi(\frac{r}{2})}{\Psi(\frac{r}{2})}=C_RN\frac{\Phi(\frac{r}{2})}{\Psi(\frac{r}{2})}.$$
Let $g=2((h-\frac{1}{2})\vee0)$, then $g\in\mathcal{F}_e$ satisfies $g(x_0)=1$, $0\le g\le 1$ in $X$, $g=0$ in $B\backslash\left(\frac{1}{2}B\right)$ and
$$\mathcal{E}(g,g)=4\mathcal{E}((h-\frac{1}{2})\vee0,(h-\frac{1}{2})\vee0)\le4\mathcal{E}(h,h).$$
By the regularity of the Dirichlet form $(\mathcal{E},\mathcal{F})$ on $L^2(X;m)$, there exists $\varphi\in\mathcal{F}\cap C_c(X)$ satisfying $0\le\varphi\le 1$ in $X$, $\varphi=1$ in an open neighborhood of $\overline{{\frac{1}{2}B}}$ and $\mathrm{supp}(\varphi)\subseteq B$. Then $\varphi g\in\mathcal{F}_e$ satisfies $(\varphi g)(x_0)=1$, $0\le\varphi g\le1$ in $X$ and $\varphi g=g1_{B}=g1_{\frac{1}{2}{B}}$. By the strongly local property and Corollary \ref{cor_attain}, we have
\begin{align*}
&R(x_0,X\backslash B)^{-1}\le\mathcal{E}(\varphi g,\varphi g)=\Gamma(\varphi g,\varphi g)(X)=\Gamma(\varphi g,\varphi g)(B)\\
&=\Gamma (g,g)(B)\le\Gamma(g,g)(X)=\mathcal{E}(g,g)\le 4\mathcal{E}(h,h)\le 4C_RN\frac{\Phi(\frac{r}{2})}{\Psi(\frac{r}{2})},
\end{align*}
hence
$$R(x_0,X\backslash B)\ge\frac{1}{4C_RN}\frac{\Psi(\frac{r}{2})}{\Phi(\frac{r}{2})}=\frac{1}{4C_RN}\frac{\frac{\Psi(\frac{r}{2})}{\Phi(\frac{r}{2})}}{\frac{\Psi(r)}{\Phi(r)}}\frac{\Psi(r)}{\Phi(r)}\ge\frac{1}{4\cdot2^{\gamma_1\vee\gamma_2}C_RN}\frac{\Psi(r)}{\Phi(r)}.$$
\end{proof}

\begin{proposition}\label{prop_main}
Let $(X,d,m,\mathcal{E},\mathcal{F})$ be an unbounded MMD space satisfying \ref{eq_VPhi},\\
\noindent \ref{eq_UHK} and \ref{eq_NLE}. Then for any ball $B=B(x_0,r)$, for any bounded $u\in\mathcal{F}_e$ which is harmonic in $B$, for any $x,y\in B$, we have
$$|u(x)-u(y)|\le\frac{R(x,y)}{R(x,X\backslash B)}\mathrm{osc}_{X\backslash B}u,$$
where
$$\mathrm{osc}_Au=\sup\nolimits_Au-\inf\nolimits_Au$$
is the oscillation of the function $u\in l(A)$ on the set $A$.
\end{proposition}

The idea of the following proof is from a recent result \cite[Theorem 6.27]{KS24a}, where a general result in the setting of $p$-energy forms was given to prove certain $p$-resistance defines indeed a metric. An alternative proof suggested by Kajino and Shimizu will be provided afterward.

\begin{proof}
Let $Y=X\backslash B$, then $Y\subsetneqq X$ is a non-empty closed subset of $X$. For fixed $x,y\in B$, we have $u-u(x)\in\mathcal{F}_e$ is harmonic in $B\backslash\{x\}$, hence $u-u(x)=h_{Y\cup\{x\}}(u-u(x))$. By the maximum principle (Corollary \ref{cor_max}), we have $u-u(x)\le(\mathrm{osc}_{Y}u)1_Y$ on $Y\cup\{x\}$, hence
$$u-u(x)=h_{Y\cup\{x\}}(u-u(x))\le h_{Y\cup\{x\}}\left((\mathrm{osc}_{Y}u)1_Y\right)=(\mathrm{osc}_{Y}u)h_{Y\cup\{x\}}(1_Y).$$
Similarly, since $u-u(x)\ge-(\mathrm{osc}_{Y}u)1_Y$ on $Y\cup\{x\}$, we have
$$u-u(x)=h_{Y\cup\{x\}}(u-u(x))\ge-h_{Y\cup\{x\}}\left((\mathrm{osc}_Yu)1_Y\right)=-\left(\mathrm{osc}_Yu\right)h_{Y\cup\{x\}}(1_Y).$$

It suffices to estimate $h_{Y\cup\{x\}}(1_Y)(y)$. Indeed, by Corollary \ref{cor_attain}, we have
\begin{align*}
&R(x,X\backslash B)^{-1}=R(x,Y)^{-1}=\mathcal{E}(h_{Y\cup\{x\}}(1_Y),h_{Y\cup\{x\}}(1_Y))\\
&=\mathcal{E}(h_{Y\cup\{x\}}(1_Y),h_{Y\cup\{x\}}(1_Y)-1)=-\mathcal{E}(h_{Y\cup\{x\}}(1_Y),h_{Y\cup\{x\}}(1_{\{x\}})).
\end{align*}
Since $\mathcal{F}_e\ni h_{Y\cup\{x,y\}}(1_{\{x\}})-h_{Y\cup\{x\}}(1_{\{x\}})=0$ on $Y\cup\{x\}$, and $h_{Y\cup\{x\}}(1_Y)\in\mathcal{F}_e$ is harmonic in $B\backslash\{x\}$, by Lemma \ref{lem_Fe} \ref{lem_Fe_har}, we have
$$\mathcal{E}(h_{Y\cup\{x\}}(1_Y),h_{Y\cup\{x,y\}}(1_{\{x\}})-h_{Y\cup\{x\}}(1_{\{x\}}))=0,$$
that is,
$$\mathcal{E}(h_{Y\cup\{x\}}(1_Y),h_{Y\cup\{x\}}(1_{\{x\}}))=\mathcal{E}(h_{Y\cup\{x\}}(1_Y),h_{Y\cup\{x,y\}}(1_{\{x\}})).$$
Let
\begin{align*}
u_1&=h_{Y\cup\{x\}}(1_Y)-h_{Y\cup\{x\}}(1_Y)(y)h_{\{x,y\}}(1_{\{y\}}),\\
v_1&=h_{Y\cup\{x,y\}}(1_{\{x\}}),
\end{align*}
then $u_1,v_1\in\mathcal{F}_e$ are harmonic in $B\backslash\{x,y\}$. Let $u_2=u_1|_{Y\cup\{x,y\}}$, $v_2=v_1|_{Y\cup\{x,y\}}$, then $u_2,v_2\in\mathcal{F}_e|_{Y\cup\{x,y\}}$ and
$$\mathcal{E}(u_1,v_1)=\mathcal{E}|_{Y\cup\{x,y\}}(u_2,v_2),$$
where $(\mathcal{E}|_{Y\cup\{x,y\}},\mathcal{F}_e|_{Y\cup\{x,y\}})$ is a resistance form on $Y\cup\{x,y\}$, which follows from Lemma \ref{lem_trace}. On $Y$, we have $u_2\in[0,1]$, $v_2=0$, $u_2\wedge v_2=0$, $u_2\vee v_2=u_2=u_2+v_2$. On $\{x\}$, we have $u_2=0$, $v_2=1$, $u_2\wedge v_2=0$, $u_2\vee v_2=1=v_2=u_2+v_2$. On $\{y\}$, we have $u_2=0$, $v_2=0$, $u_2\wedge v_2=0$, $u_2\vee v_2=0=u_2+v_2$. Hence $u_2\wedge v_2=0$, $u_2\vee v_2=u_2+v_2$. By the following easy result, Lemma \ref{lem_ele}, we have
$$\mathcal{E}(u_1,v_1)=\mathcal{E}|_{Y\cup\{x,y\}}(u_2,v_2)\le0,$$
that is,
$$\mathcal{E}(h_{Y\cup\{x\}}(1_Y), h_{Y\cup\{x,y\}}(1_{\{x\}}))\le\mathcal{E}(h_{Y\cup\{x\}}(1_Y)(y)h_{\{x,y\}}(1_{\{y\}}), h_{Y\cup\{x,y\}}(1_{\{x\}})),$$
hence
\begin{align*}
&R(x,X\backslash B)^{-1}=-\mathcal{E}(h_{Y\cup\{x\}}(1_Y),h_{Y\cup\{x\}}(1_{\{x\}}))=-\mathcal{E}(h_{Y\cup\{x\}}(1_Y),h_{Y\cup\{x,y\}}(1_{\{x\}}))\\
&\ge-\mathcal{E}(h_{Y\cup\{x\}}(1_Y)(y)h_{\{x,y\}}(1_{\{y\}}), h_{Y\cup\{x,y\}}(1_{\{x\}}))\\
&=-h_{Y\cup\{x\}}(1_Y)(y)\mathcal{E}(h_{\{x,y\}}(1_{\{y\}}), h_{Y\cup\{x,y\}}(1_{\{x\}})).
\end{align*}
Since $\mathcal{F}_e\ni h_{Y\cup\{x,y\}}(1_{\{x\}})-h_{\{x,y\}}(1_{\{x\}})=0$ on $\{x,y\}$ and $h_{\{x,y\}}(1_{\{y\}})\in\mathcal{F}_e$ is harmonic in $X\backslash\{x,y\}$, by Lemma \ref{lem_Fe} \ref{lem_Fe_har}, we have
$$\mathcal{E}(h_{\{x,y\}}(1_{\{y\}}),h_{Y\cup\{x,y\}}(1_{\{x\}})-h_{\{x,y\}}(1_{\{x\}}))=0,$$
that is,
$$\mathcal{E}(h_{\{x,y\}}(1_{\{y\}}), h_{Y\cup\{x,y\}}(1_{\{x\}}))=\mathcal{E}(h_{\{x,y\}}(1_{\{y\}}), h_{\{x,y\}}(1_{\{x\}})).$$
Hence
\begin{align*}
&R(x,X\backslash B)^{-1}\ge-h_{Y\cup\{x\}}(1_Y)(y)\mathcal{E}(h_{\{x,y\}}(1_{\{y\}}), h_{Y\cup\{x,y\}}(1_{\{x\}}))\\
&=-h_{Y\cup\{x\}}(1_Y)(y)\mathcal{E}(h_{\{x,y\}}(1_{\{y\}}), h_{\{x,y\}}(1_{\{x\}}))\\
&=h_{Y\cup\{x\}}(1_Y)(y)\mathcal{E}(1-h_{\{x,y\}}(1_{\{y\}}), h_{\{x,y\}}(1_{\{x\}}))\\
&=h_{Y\cup\{x\}}(1_Y)(y)\mathcal{E}(h_{\{x,y\}}(1_{\{x\}}), h_{\{x,y\}}(1_{\{x\}}))\\
&=h_{Y\cup\{x\}}(1_Y)(y)R(x,y)^{-1},
\end{align*}
which implies
$$h_{Y\cup\{x\}}(1_Y)(y)\le\frac{R(x,y)}{R(x,X\backslash B)}.$$

Therefore, we have
$$|u(x)-u(y)|\le\frac{R(x,y)}{R(x,X\backslash B)}\mathrm{osc}_Yu.$$
\end{proof}

\begin{lemma}\label{lem_ele}
Let $X$ be a set. Let $(\mathcal{E},\mathcal{G})$ be a resistance form on $X$. If $u,v\in\mathcal{G}$ satisfy $u\wedge v=0$ and $u\vee v=u+v$, then $\mathcal{E}(u,v)\le0$.
\end{lemma}

\begin{proof}
By the parallelogram identity and \ref{def_RF5}, we have
$$\mathcal{E}(u\wedge v,u\wedge v)+\mathcal{E}(u\vee v,u\vee v)\le\mathcal{E}(u,u)+\mathcal{E}(v,v).$$
By assumption, we have
$$\mathcal{E}(u,u)+\mathcal{E}(v,v)\ge\mathcal{E}(u+v,u+v)=\mathcal{E}(u,u)+\mathcal{E}(v,v)+2\mathcal{E}(u,v),$$
hence $\mathcal{E}(u,v)\le0$.
\end{proof}

The alternative proof relies on the following Lipschitz continuity result for Green functions in the setting of Kigami's resistance form theory.

\begin{lemma}[{\cite[THEOREM 4.1]{Kig12}}]\label{lem_Green}
Let $(X,d,m,\mathcal{E},\mathcal{F})$ be an unbounded MMD space satisfying \ref{eq_VPhi}, \ref{eq_UHK} and \ref{eq_NLE}. Let $D\subsetneqq X$ be a non-empty open subset of $X$, and $\mathcal{F}_{e,D}$ given by Equation (\ref{eq_FeD}). Then $(\mathcal{E},\mathcal{F}_{e,D})$ is a Hilbert space and there exists a unique $g_D:X\times X\to \mathbb{R}$ satisfying that for any $x\in X$, $g_D(x,\cdot)\in \mathcal{F}_{e,D}$ and $\mathcal{E}(g_D(x,\cdot),u)=u(x)$ for any $u\in \mathcal{F}_{e,D}$. Moreover, for any $x,y,z\in X$, we have $g_D(x,x)=R(x,X\backslash D)$ and
$$\lvert g_D(x,y)-g_D(x,z)\rvert\le R(y,z).$$
\end{lemma}

\begin{proof}[An alternative proof of Proposition \ref{prop_main} using Lemma \ref{lem_Green}]
Let $Y=X\backslash B$ and $x\in B$ as in the previous proof. Recall that it suffices to estimate $h_{Y\cup\{x\}}(1_Y)$. Indeed, let $h=h_{Y\cup\{x\}}(1_{\{x\}})=1-h_{Y\cup\{x\}}(1_Y)$, then $h\in \mathcal{F}_{e,B}$. By Corollary \ref{cor_attain}, we have $\mathcal{E}(h,h)=R(x,X\backslash B)^{-1}$. For any $v\in \mathcal{F}_{e,B}$, we have $v-v(x)h\in \mathcal{F}_{e,B\backslash\{x\}}$. Since $h$ is harmonic in $B\backslash\{x\}$, by Lemma \ref{lem_Fe} \ref{lem_Fe_har}, we have $\mathcal{E}(h,v-v(x)h)=0$, that is,
$$\mathcal{E}(h,v)=v(x)\mathcal{E}(h,h)=\frac{v(x)}{R(x,X\backslash B)}\text{ for any }v\in \mathcal{F}_{e,B}.$$
By Lemma \ref{lem_Green}, we have $h=\frac{g_B(x,\cdot)}{R(x,X\backslash B)}$, $g_B(x,x)=R(x,X\backslash B)$, and
\begin{align*}
h_{Y\cup\{x\}}(1_Y)=1-h=\frac{g_B(x,x)-g_B(x,\cdot)}{R(x,X\backslash B)}\le \frac{R(x,\cdot)}{R(x,X\backslash B)}.
\end{align*}
The remainder of the argument proceeds as in the previous proof.
\end{proof}

We need the following $L^1$-mean value inequality.

\begin{lemma}[{\cite[THEOREM 6.3, LEMMA 9.2]{GHL15}}]\label{lem_MV}
Let $(X,d,m,\mathcal{E},\mathcal{F})$ be an unbounded MMD space satisfying \ref{eq_VD}, \ref{eq_FK} and \ref{eq_CS}. Then there exists $C\in(0,+\infty)$ such that for any ball $B=B(x_0,r)$, for any $u\in\mathcal{F}$ which is harmonic in $2B$, we have
$$\lVert {u}\rVert_{L^\infty(B)}\le C\dashint_{2B}|u|\mathrm{d}m.$$
\end{lemma}

We now give the proof of ``\ref{eq_NLE}$\Rightarrow$\ref{eq_HR}" as follows.

\begin{proof}[Proof of ``\ref{eq_NLE}$\Rightarrow$\ref{eq_HR}"]
Note that we need to prove \ref{eq_HR} from \ref{eq_UHK} and \ref{eq_NLE}. By Proposition \ref{prop_res}, we have \ref{eq_res}. Let $B=B(x_0,r)$ and $u\in\mathcal{F}$ harmonic in $2B$. By \ref{eq_VD} and Lemma \ref{lem_MV}, we have
$$M:=\lVert {u}\rVert_{L^\infty(\frac{3}{2}B)}\lesssim\dashint_{2B}|u|\mathrm{d}m.$$
Let $v=(u\vee(-M))\wedge M$, then $v\in\mathcal{F}\cap L^\infty(X;m)$, $\lVert {v}\rVert_{L^\infty(X)}\le M$ and $v=u$ in $\frac{3}{2}B$, hence $v\in\mathcal{F}\subseteq\mathcal{F}_e$ is bounded in $X$ and harmonic in $\frac{3}{2}B$. By Proposition \ref{prop_main}, for any $x,y\in B\subseteq\frac{3}{2}B$, we have
\begin{align*}
&|v(x)-v(y)|\le\frac{R(x,y)}{R(x,X\backslash\left(\frac{3}{2}B\right))}\mathrm{osc}_{X\backslash\left(\frac{3}{2}B\right)}v\\
&\le 2\frac{R(x,y)}{R(x,X\backslash B(x,\frac{r}{2}))}\lVert {v}\rVert_{L^\infty(X)}\le2M\frac{R(x,y)}{R(x,X\backslash B(x,\frac{r}{2}))},
\end{align*}
hence
$$|u(x)-u(y)|\le2 \lVert {u}\rVert_{L^\infty(\frac{3}{2}B)}\frac{R(x,y)}{R(x,X\backslash B(x,\frac{r}{2}))}\lesssim \frac{R(x,y)}{R(x,X\backslash B(x,\frac{r}{2}))}\dashint_{2B}|u|\mathrm{d}m.$$
By \ref{eq_res} and Lemma \ref{lem_ptball}, we have
$$|u(x)-u(y)|\lesssim\frac{\left(\frac{\Psi}{\Phi}\right)(d(x,y))}{\left(\frac{\Psi}{\Phi}\right)(\frac{r}{2})}\dashint_{2B}|u|\mathrm{d}m\asymp\frac{\left(\frac{\Psi}{\Phi}\right)(d(x,y))}{\left(\frac{\Psi}{\Phi}\right)(r)}\dashint_{2B}|u|\mathrm{d}m,$$
hence \ref{eq_HR} holds.
\end{proof}

\begin{remark}
A natural question concerns the sharpness of the H\"older exponents $\beta_1-\alpha_1$ and $\beta_2-\alpha_2$ arising from \ref{eq_HR}. Indeed, on a class of unbounded p.c.f. self-similar sets and the corresponding cable systems, these exponents are sharp, see the proof of \cite[Proposition 4.1, Proposition 4.3]{DRY23} for the Vicsek and the Sierpi\'nski gasket cable systems. However, on the unbounded Sierpi\'nski carpet and the Sierpi\'nski carpet cable system, it was conjectured in \cite[Conjecture 5.4, Open Question 5.5]{ABCRST3} that these exponents are \emph{not} sharp, and a new fractal dimension, the so-called topological Hausdorff dimension (see \cite{BBE15}), should be involved.
\end{remark}

\section{{\normalsize{Some examples related to the Sierpi\'nski carpet and the Vicsek set}}}\label{sec_eg}

In this section, we give some examples related to the Sierpi\'nski carpet and the Vicsek set, including blowups of one fractal by another. In $\mathbb{R}^2$, let
\begin{equation}\label{eq_SCps}
\begin{aligned}
&p_1=(0,0),p_2=(\frac{1}{2},0),p_3=(1,0),p_4=(1,\frac{1}{2}),\\
&p_5=(1,1),p_6=(\frac{1}{2},1),p_7=(0,1),p_8=(0,\frac{1}{2}),
\end{aligned}
\end{equation}
and
$$g_i(x)=\frac{1}{3}\left(x-p_i\right)+p_i,x\in \mathbb{R}^2,i=1,\ldots,8.$$
Then the Sierpi\'nski carpet is the unique non-empty compact set $K$ in $\mathbb{R}^2$ satisfying $K=\cup_{i=1}^8g_i(K)$, see Figure \ref{fig_SC}. Let $\alpha=\log8/\log3$. By \cite[Theorem 1.5.7, Proposition 1.5.8]{Kig01book}, the Hausdorff dimension of $K$ is $\alpha$ and the $\alpha$-dimensional Hausdorff measure $\mathcal{H}^\alpha(K)\in(0,+\infty)$. The construction of Brownian motion on the Sierpi\'nski carpet was given by \cite{BB89} and \cite{KZ92} using approximations of stochastic processes. See also \cite{GY19} for another construction using $\Gamma$-convergence of quadratic forms, \cite{BB99} for the construction on higher dimensional Sierpi\'nski carpets, \cite{CQ24} for the construction on unconstrained Sierpi\'nski carpets. Let $\beta=\log(8\rho)/\log3\approx2.09697$, where $\rho\approx1.25147\in[\frac{7}{6},\frac{3}{2}]$ is some parameter from resistance estimates in \cite{BB90,BBS90}.

For $a\in \mathbb{R}$, $A\subseteq \mathbb{R}^2$, denote $aA=\{ax:x\in A\}$.

Firstly, we consider the unbounded Sierpi\'nski carpet $K^\infty=\cup_{n\ge0}3^nK$ as follows, where  $\{3^nK\}_{n\ge0}$ is an increasing sequence of subsets of $\mathbb{R}^2$. Let $\alpha_1=\alpha_2=\alpha$, $\beta_1=\beta_2=\beta$. It is obvious that \ref{eq_VPhi} holds. It was proved in \cite[Theorem 1.1]{BB92} that \hyperlink{eq_HKPsi}{HK($\Psi$)} and \ref{eq_HHK} hold. In a general setting of fractional diffusions introduced in \cite[Section 3]{Bar98}, it was proved in \cite[Theorem 3.7, Lemma 3.4]{ABCRST3} that \ref{eq_wBE} and \ref{eq_HHK} hold. It was proved in \cite[Theorem 8.3]{BB92} that all functions in the domain of the generator of the Dirichlet form are $(\beta-\alpha)$-H\"older continuous. By the main result, Theorem \ref{thm_main}, we now have the H\"older regularity for harmonic functions \ref{eq_HR} and the H\"older estimate \ref{eq_HHKexp} with exponential terms for the heat kernel.

\begin{figure}[ht]
\centering
\begin{subfigure}[b]{0.2\textwidth}
\centering
\begin{tikzpicture}[scale=0.15]
\draw(0,0)--(2,0)--(2,2)--(0,2)--cycle;

\draw[fill=black] (0,0) circle (0.25);
\draw[fill=black] (0,2) circle (0.25);
\draw[fill=black] (2,2) circle (0.25);
\draw[fill=black] (2,0) circle (0.25);

\draw[fill=black] (1,0) circle (0.25);
\draw[fill=black] (2,1) circle (0.25);
\draw[fill=black] (1,2) circle (0.25);
\draw[fill=black] (0,1) circle (0.25);
\end{tikzpicture}
\caption{$V^{(0)}$}
\end{subfigure}
\hspace{2em}
\begin{subfigure}[b]{0.2\textwidth}
\centering
\begin{tikzpicture}[scale=0.15]

\foreach\x in {0}{
\foreach\y in {0}{

\foreach\xx in {0,2,4}{
\foreach\yy in {0,4}{
\draw(\x+\xx+0,\y+\yy+0)--(\x+\xx+2,\y+\yy+0)--(\x+\xx+2,\y+\yy+2)--(\x+\xx+0,\y+\yy+2)--cycle;

\draw[fill=black] (\x+\xx+0,\y+\yy+0) circle (0.25);
\draw[fill=black] (\x+\xx+0,\y+\yy+2) circle (0.25);
\draw[fill=black] (\x+\xx+2,\y+\yy+2) circle (0.25);
\draw[fill=black] (\x+\xx+2,\y+\yy+0) circle (0.25);

\draw[fill=black] (\x+\xx+1,\y+\yy+0) circle (0.25);
\draw[fill=black] (\x+\xx+2,\y+\yy+1) circle (0.25);
\draw[fill=black] (\x+\xx+1,\y+\yy+2) circle (0.25);
\draw[fill=black] (\x+\xx+0,\y+\yy+1) circle (0.25);
}
}

\foreach\xx in {0,4}{
\foreach\yy in {2}{
\draw(\x+\xx+0,\y+\yy+0)--(\x+\xx+2,\y+\yy+0)--(\x+\xx+2,\y+\yy+2)--(\x+\xx+0,\y+\yy+2)--cycle;

\draw[fill=black] (\x+\xx+0,\y+\yy+0) circle (0.25);
\draw[fill=black] (\x+\xx+0,\y+\yy+2) circle (0.25);
\draw[fill=black] (\x+\xx+2,\y+\yy+2) circle (0.25);
\draw[fill=black] (\x+\xx+2,\y+\yy+0) circle (0.25);

\draw[fill=black] (\x+\xx+1,\y+\yy+0) circle (0.25);
\draw[fill=black] (\x+\xx+2,\y+\yy+1) circle (0.25);
\draw[fill=black] (\x+\xx+1,\y+\yy+2) circle (0.25);
\draw[fill=black] (\x+\xx+0,\y+\yy+1) circle (0.25);
}
}

}
}
\end{tikzpicture}
\caption{$V^{(1)}$}
\end{subfigure}
\hspace{2em}
\begin{subfigure}[b]{0.3\textwidth}
\centering
\begin{tikzpicture}[scale=0.15]

\foreach\x in {0,6,12}{
\foreach\y in {0,12}{

\foreach\xx in {0,2,4}{
\foreach\yy in {0,4}{
\draw(\x+\xx+0,\y+\yy+0)--(\x+\xx+2,\y+\yy+0)--(\x+\xx+2,\y+\yy+2)--(\x+\xx+0,\y+\yy+2)--cycle;

\draw[fill=black] (\x+\xx+0,\y+\yy+0) circle (0.25);
\draw[fill=black] (\x+\xx+0,\y+\yy+2) circle (0.25);
\draw[fill=black] (\x+\xx+2,\y+\yy+2) circle (0.25);
\draw[fill=black] (\x+\xx+2,\y+\yy+0) circle (0.25);

\draw[fill=black] (\x+\xx+1,\y+\yy+0) circle (0.25);
\draw[fill=black] (\x+\xx+2,\y+\yy+1) circle (0.25);
\draw[fill=black] (\x+\xx+1,\y+\yy+2) circle (0.25);
\draw[fill=black] (\x+\xx+0,\y+\yy+1) circle (0.25);
}
}

\foreach\xx in {0,4}{
\foreach\yy in {2}{
\draw(\x+\xx+0,\y+\yy+0)--(\x+\xx+2,\y+\yy+0)--(\x+\xx+2,\y+\yy+2)--(\x+\xx+0,\y+\yy+2)--cycle;

\draw[fill=black] (\x+\xx+0,\y+\yy+0) circle (0.25);
\draw[fill=black] (\x+\xx+0,\y+\yy+2) circle (0.25);
\draw[fill=black] (\x+\xx+2,\y+\yy+2) circle (0.25);
\draw[fill=black] (\x+\xx+2,\y+\yy+0) circle (0.25);

\draw[fill=black] (\x+\xx+1,\y+\yy+0) circle (0.25);
\draw[fill=black] (\x+\xx+2,\y+\yy+1) circle (0.25);
\draw[fill=black] (\x+\xx+1,\y+\yy+2) circle (0.25);
\draw[fill=black] (\x+\xx+0,\y+\yy+1) circle (0.25);
}
}
}
}

\foreach\x in {0,12}{
\foreach\y in {6}{

\foreach\xx in {0,2,4}{
\foreach\yy in {0,4}{
\draw(\x+\xx+0,\y+\yy+0)--(\x+\xx+2,\y+\yy+0)--(\x+\xx+2,\y+\yy+2)--(\x+\xx+0,\y+\yy+2)--cycle;

\draw[fill=black] (\x+\xx+0,\y+\yy+0) circle (0.25);
\draw[fill=black] (\x+\xx+0,\y+\yy+2) circle (0.25);
\draw[fill=black] (\x+\xx+2,\y+\yy+2) circle (0.25);
\draw[fill=black] (\x+\xx+2,\y+\yy+0) circle (0.25);

\draw[fill=black] (\x+\xx+1,\y+\yy+0) circle (0.25);
\draw[fill=black] (\x+\xx+2,\y+\yy+1) circle (0.25);
\draw[fill=black] (\x+\xx+1,\y+\yy+2) circle (0.25);
\draw[fill=black] (\x+\xx+0,\y+\yy+1) circle (0.25);
}
}

\foreach\xx in {0,4}{
\foreach\yy in {2}{
\draw(\x+\xx+0,\y+\yy+0)--(\x+\xx+2,\y+\yy+0)--(\x+\xx+2,\y+\yy+2)--(\x+\xx+0,\y+\yy+2)--cycle;

\draw[fill=black] (\x+\xx+0,\y+\yy+0) circle (0.25);
\draw[fill=black] (\x+\xx+0,\y+\yy+2) circle (0.25);
\draw[fill=black] (\x+\xx+2,\y+\yy+2) circle (0.25);
\draw[fill=black] (\x+\xx+2,\y+\yy+0) circle (0.25);

\draw[fill=black] (\x+\xx+1,\y+\yy+0) circle (0.25);
\draw[fill=black] (\x+\xx+2,\y+\yy+1) circle (0.25);
\draw[fill=black] (\x+\xx+1,\y+\yy+2) circle (0.25);
\draw[fill=black] (\x+\xx+0,\y+\yy+1) circle (0.25);
}
}
}
}

\end{tikzpicture}
\caption{$V^{(2)}$}
\end{subfigure}
\caption{$V^{(0)}$, $V^{(1)}$ and $V^{(2)}$}\label{fig_V012}
\end{figure}
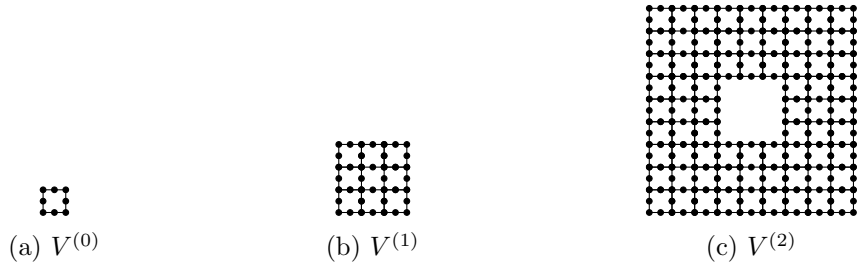

Secondly, we consider the Sierpi\'nski carpet cable system. Let $V_0=\{p_1,\ldots,p_8\}$ and $V_{n+1}=\cup_{i=1}^8g_i(V_n)$ for any $n\ge0$, then $\{V_n\}_{n\ge0}$ is an increasing sequence of finite subsets of $K$, and the closure of $\cup_{n\ge0}V_n$ with respect to the Euclidean topology is $K$. For any $n\ge0$, let $V^{(n)}=3^nV_n$, see Figure \ref{fig_V012} for $V^{(0)}$, $V^{(1)}$ and $V^{(2)}$.

Obviously, $\{V^{(n)}\}_{n\ge0}$ is an increasing sequence of finite sets. Let $V=\cup_{n\ge0}V^{(n)}$ and
$$E=\left\{\{p,q\}:p,q\in V,|p-q|=\frac{1}{2}\right\},$$
then $(V,E)$ is an infinite, locally bounded, connected (undirected) graph, whose corresponding cable system is the Sierpi\'nski carpet cable system, see \cite[Section 3]{DRY23} for a more detailed construction. Let $\alpha_1=1$, $\beta_1=2$, $\alpha_2=\alpha$, $\beta_2=\beta$. It is obvious that \ref{eq_VPhi} holds. Using the stability result \cite{BBK06} for heat kernel estimates or parabolic Harnack inequalities, we have \hyperlink{eq_HKPsi}{HK($\Psi$)} also holds, which is similar to the case on the unbounded Sierpi\'nski carpet. By the main result, Theorem \ref{thm_main}, we have \ref{eq_HHK}, \ref{eq_HHKexp}, \ref{eq_wBE} and \ref{eq_HR} all hold. In particular, the validity of \ref{eq_HR} implies the validity of \ref{eq_GRH}, which was left open in \cite{DRY23}. Moreover, we also have \ref{eq_GHK}, which can also be given directly by \ref{eq_HHKexp}.

Thirdly, we consider blowups of fractals related to the Sierpi\'nski carpet and the Vicsek set. However, we blowup one fractal by a different fractal. In what follows, we work in $\mathbb{R}^2$.

Let $K,M$ be one of the following fractals in $\mathbb{R}^2$: the unit interval $[0,1]\times\{0\}$, the Vicsek set, and the Sierpi\'nski carpet. Let $\alpha_1$ and $\beta_1$ be the Hausdorff dimension and the walk dimension of $K$, and let $\alpha_2$ and $\beta_2$ be the Hausdorff dimension and the walk dimension of $M$. Here, $K$ serves as the ``cell" fractal and $M$ provides the ``model" to blowup. Let $(\mathcal{E}^K,\mathcal{F}^K)$ be the strongly local regular Dirichlet form on $L^2(K;\mu)$ with two-sided sub-Gaussian heat kernel estimates with walk dimension $\beta_1$, where $\mu$ is the normalized $\alpha_1$-Hausdorff measure on $K$, see \cite{Lindstrom90, Kig01book, Str06book} for the case $K$ is the Vicsek set.

Let $\{g^{(1)}_i\}_{i=1,\ldots,N^{(1)}}$, $\{g^{(2)}_i\}_{i=1,\ldots,N^{(2)}}$ be the collections of contraction mappings generating $K$, $M$, respectively, that is, for $k=1,2$, $N^{(k)},L^{(k)}$ are two positive integers, $p^{(k)}_1=(0,0)$, \ldots, $p^{(k)}_{N^{(k)}}\in \mathbb{R}^2$,
$$g^{(k)}_i(x)=\frac{1}{L^{(k)}}x+\frac{L^{(k)}-1}{L^{(k)}}p_i,x\in \mathbb{R}^2,i=1,\ldots, N^{(k)},$$
and $K$, $M$ are the unique non-empty compact sets in $\mathbb{R}^2$ satisfying $K=\cup_{i=1}^{N^{(1)}}g^{(1)}_i(K)$, $M=\cup_{i=1}^{N^{(2)}}g^{(2)}_i(M)$, respectively.

\begin{itemize}
\item For the unit interval, we have $N^{(k)}=2$, $L^{(k)}=2$ and
$$p^{(k)}_1=(0,0), p^{(k)}_2=(1,0).$$
\item For the Vicsek set, we have $N^{(k)}=5$, $L^{(k)}=3$ and
$$p^{(k)}_1=(0,0),p^{(k)}_2=(1,0),p^{(k)}_3=(1,1),p^{(k)}_4=(0,1),p^{(k)}_5=\left(\frac{1}{2},\frac{1}{2}\right).$$
\item For the Sierpi\'nski carpet, we have $N^{(k)}=8$, $L^{(k)}=3$ and $p^{(k)}_1,\ldots,p^{(k)}_8$ are $p_1,\ldots,p_8$ from Equation (\ref{eq_SCps}).
\end{itemize}

We assume that for $i,j=1,\ldots,N^{(2)}$, if $g^{(2)}_i\left([0,1]^2\right)\cap g^{(2)}_j\left([0,1]^2\right)\ne\emptyset$, then $g^{(2)}_i(K)\cap g^{(2)}_j(K)\ne\emptyset$. Let $K_0=K$. For any $n\ge0$, let $K_{n+1}=\cup_{i=1}^{N^{(2)}}g^{(2)}_i(K_n)$, $K^{(n)}=(L^{(2)})^nK_n$, then
\begin{align*}
&K^{(n+1)}=(L^{(2)})^{n+1}K_{n+1}=(L^{(2)})^{n+1}\left(\cup_{i=1}^{N^{(2)}}g^{(2)}_i(K_n)\right)\\
&=(L^{(2)})^{n+1}\left(g^{(2)}_1(K_n)\cup\bigcup_{i=2}^{N^{(2)}}g^{(2)}_i(K_n)\right)\\
&=\left((L^{(2)})^nK_n\right)\cup\left((L^{(2)})^{n+1}\left(\bigcup_{i=2}^{N^{(2)}}g^{(2)}_i(K_n)\right)\right)\\
&\supseteq (L^{(2)})^nK_n=K^{(n)},
\end{align*}
hence $\{K^{(n)}\}_{n\ge0}$ is an increasing sequence of subsets of $\mathbb{R}^2$.

Let $X=\cup_{n\ge0}K^{(n)}$. We say that $W\subseteq X$ is a cell if $W$ is a translation of $K^{(0)}=K_0=K$. Let $d$ be the Euclidean metric in $\mathbb{R}^2$. Let $m$ be the measure on $X$ satisfying that for any cell $W$, $m|_{W}$ is the normalized $\alpha_1$-Hausdorff measure on $W$. For any cell $W$, by translation, there exists a strongly local regular Dirichlet form $(\mathcal{E}^W,\mathcal{F}^W)$ on $L^2(W;m|_W)$. Let
$$\mathcal{K}=\left\{u\in C_c(X):u|_W\in\mathcal{F}^W\text{ for any cell }W\right\}.$$
Define
\begin{align*}
\mathcal{E}(u,u)&=\sum_{W\text{ is a cell}}\mathcal{E}^W(u|_W,u|_W),\\
\mathcal{F}&=\mathcal{E}_1\text{-closure of }\mathcal{K},
\end{align*}
then $(\mathcal{E},\mathcal{F})$ is a strongly local regular Dirichlet form on $L^2(X;m)$. Obviously, $(X,d,m,\mathcal{E},\mathcal{F})$ is an unbounded MMD space satisfying \ref{eq_VPhi}. Since two-sided heat kernel estimates hold on $K,M$, by \cite[Theorem 3.14]{GH14a}, we have the elliptic Harnack inequalities and the resistance estimates both hold on $K,M$. Combining these two conditions on $K,M$, we obtain easily the corresponding elliptic Harnack inequalities and the resistance estimates on $X$. Then by \cite[Theorem 3.14]{GH14a} again, we have \hyperlink{eq_HKPsi}{HK($\Psi$)} holds on $X$. Indeed, in the strongly recurrent setting, only resistance estimates are enough for the above argument, see \cite{BCK05, Kum04, Hu08}. By the main result, Theorem \ref{thm_main}, we have \ref{eq_HHK}, \ref{eq_HHKexp}, \ref{eq_wBE} and \ref{eq_HR} all hold.

We give some examples as follows.

\begin{enumerate}[label=(\arabic*)]
\item Let $K$ be the Sierpi\'nski carpet and $M$ the unit interval, see Figure \ref{fig_SC01}.

\begin{figure}[ht]
\begin{center}
\begin{tikzpicture}
\foreach\x in {50pt}{
\foreach\y in {0,1,2,3,4,5}{
\node at (\y*\x,0) {\includegraphics[width=\x]{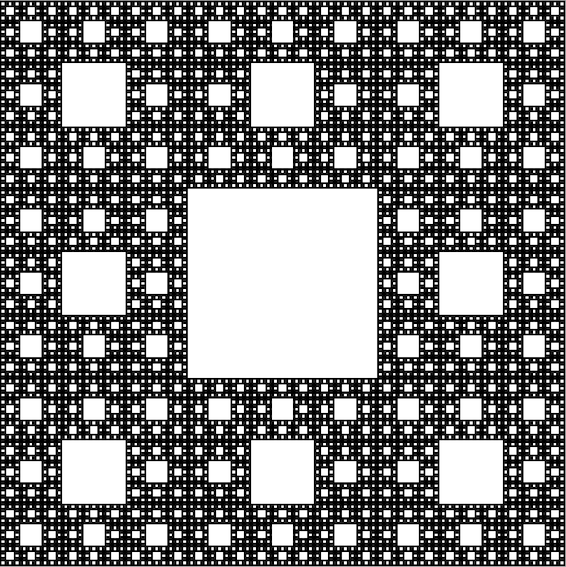}};
}
}
\draw[fill=black] (290pt,0pt) circle (1pt);
\draw[fill=black] (295pt,0pt) circle (1pt);
\draw[fill=black] (300pt,0pt) circle (1pt);
\end{tikzpicture}
\end{center}
\caption{Blow up the Sierpi\'nski carpet by the unit interval}\label{fig_SC01}
\end{figure}
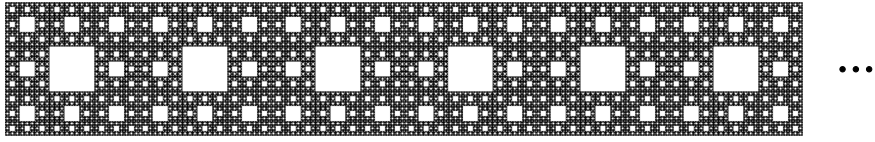

\item Let $K$ be the Sierpi\'nski carpet and $M$ the Vicsek set, see Figure \ref{fig_SCVC} for $K^{(2)}$.

\begin{figure}[ht]
\begin{center}
\begin{tikzpicture}
\foreach\w in {20pt}{
\foreach\x in {0,1,2}{
\foreach\y in {0,1,2}{
\pgfmathparse{int(mod(\x+\y,2))}
\let\r\pgfmathresult
\ifnum\r=0
\foreach\xx in {0,1,2}{
\foreach\yy in {0,1,2}{
\pgfmathparse{int(mod(\xx+\yy,2))}
\let\rr\pgfmathresult
\ifnum\rr=0
\node at (\x*3*\w+\xx*\w,\y*3*\w+\yy*\w) {\includegraphics[width=\w]{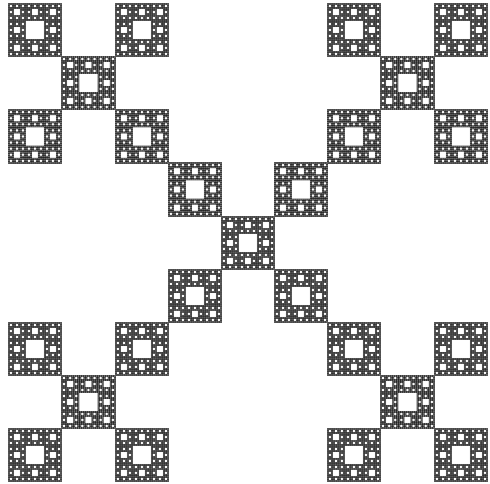}};
\fi
}
}
\fi
}
}
}
\end{tikzpicture}
\end{center}
\caption{Blow up the Sierpi\'nski carpet by the Vicsek set}\label{fig_SCVC}
\end{figure}

\item Let $K$ be the Vicsek set and $M$ the Sierpi\'nski carpet, see Figure \ref{fig_VCSC} for $K^{(2)}$.

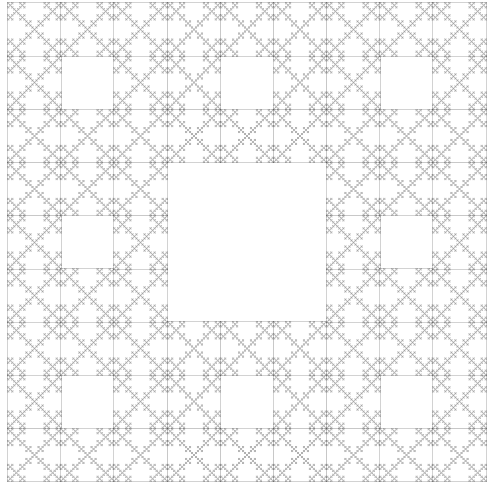
\begin{figure}[htp]
\begin{center}
\begin{tikzpicture}
\foreach\w in {20pt}{
\foreach\x in {0,1,2}{
\foreach\y in {0,1,2}{
\pgfmathparse{\x==1 && \y ==1 ? int(1) : int(0)}
\ifnum\pgfmathresult=1
\else
{
\foreach\xx in {0,1,2}{
\foreach\yy in {0,1,2}{
\pgfmathparse{\xx==1 && \yy ==1 ? int(1) : int(0)}
\ifnum\pgfmathresult=1
\else
\node at (\x*3*\w+\xx*\w,\y*3*\w+\yy*\w) {\includegraphics[width=20.2pt]{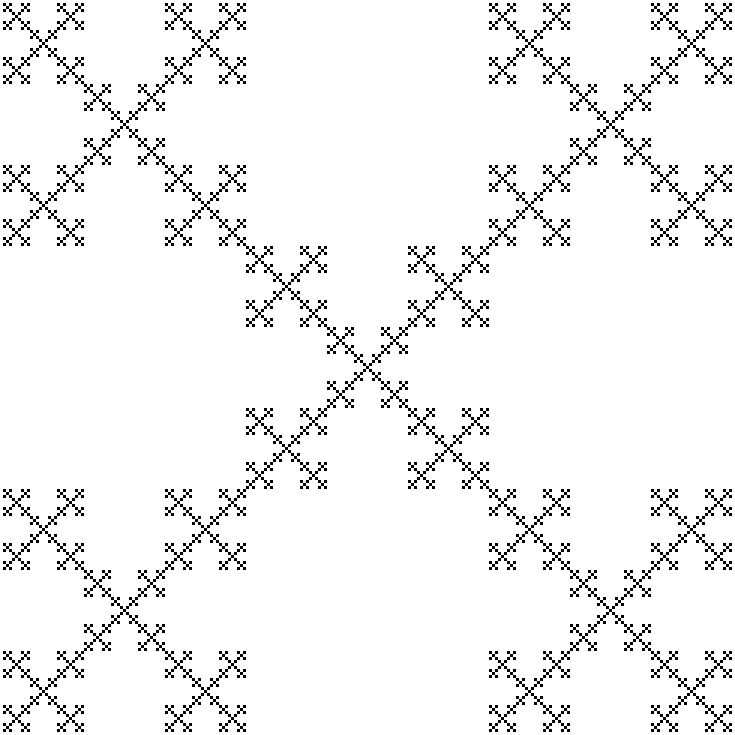}};
\fi
}
}
}
\fi
}
}
}
\end{tikzpicture}
\end{center}
\caption{Blow up the Vicsek set by the Sierpi\'nski carpet}\label{fig_VCSC}
\end{figure}

Since $(\Psi/\Phi)(r)=r$ for $r\in(0,1)$, \ref{eq_HHKexp} gives
\begin{align}
&|p_t(x,y_1)-p_t(x,y_2)|\nonumber\\
&\le C_1\frac{d(y_1,y_2)}{t}\left(\exp\left(-\Upsilon\left(C_2d(x,y_1),t\right)\right)+\exp\left(-\Upsilon\left(C_2d(x,y_2),t\right)\right)\right)\label{eq_vcsc_lip}
\end{align}
for any $x,y_1,y_2\in X$ with $d(y_1,y_2)<1$.

Let $S_0=[(0,0),(1,1)]\cup[(0,1),(1,0)]$ be the union of the diagonals of the unit square $[0,1]^2\subseteq \mathbb{R}^2$. For any $n\ge0$, let $S_{n+1}=\cup_{i=1}^{N^{(1)}}g_i^{(1)}(S_n)$. Then $\{S_n\}_{n\ge0}$ is an increasing sequence of subsets of $K$, and the closure of $\cup_{n\ge0}S_n$ is $K$. Let $\mathcal{S}_0=\cup_{n\ge0}S_n$, $\nu^{(0)}$ the measure on $\mathcal{S}_0$ such that for any $n\ge0$, $\nu^{(0)}|_{S_n}$ is the normalized 1-dimensional Hausdorff measure satisfying
$$\nu^{(0)}\left([(0,0),(1,1)]\right)=\nu^{(0)}\left([(0,1),(1,0)]\right)=1.$$
For any $n\ge0$, let $\mathcal{S}_{n+1}=\cup_{i=1}^{N^{(2)}}g_i^{(2)}(\mathcal{S}_n)$, $\mathcal{S}^{(n)}=(L^{(2)})^n\mathcal{S}_n$, then similar to $\{K^{(n)}\}_{n\ge0}$, $\{\mathcal{S}^{(n)}\}_{n\ge0}$ is an increasing sequence of subsets of $\mathbb{R}^2$, let $\mathcal{S}=\cup_{n\ge0}\mathcal{S}^{(n)}$. For any cell $W\subseteq X$, we have $\mathcal{S}_W:=\mathcal{S}\cap W$ is a translation of $\mathcal{S}^{(0)}=\mathcal{S}_0$. For any two cells $W_1,W_2\subseteq X$, we have either the cardinal $\#(\mathcal{S}_{W_1}\cap\mathcal{S}_{W_2})\in\{0,1\}$ or $\mathcal{S}_{W_1}\cap\mathcal{S}_{W_2}$ is countably infinite, hence there exists a unique $\nu$ on $\mathcal{S}$ such that for any cell $W$, $\nu|_{\mathcal{S}_W}$ is the push forward of $\nu^{(0)}$ under the translation from $\mathcal{S}^{(0)}=\mathcal{S}_0$ to $\mathcal{S}_W$. Note that $\nu(B(x_0),r)=+\infty$ for any ball $B(x_0,r)$ and $\nu\perp m$, that is, $\nu$ and $m$ are mutually singular.

By the Lipschitz estimate, Equation (\ref{eq_vcsc_lip}), we have the weak gradient $\partial p_t(x,\cdot)$ exists $\nu$-a.e. for any $t\in(0,+\infty)$, for any $x\in X$, and has the following estimate. For any $t\in(0,+\infty)$, for any $x\in X$, for $\nu$-a.e. $y\in\mathcal{S}$, we have
\begin{align}
&|\partial_yp_t(x,y)|\le\frac{C_1}{t}\exp\left(-\Upsilon\left(C_2d(x,y),t\right)\right)\nonumber\\
&\le
\begin{cases}
\frac{C_3}{t}\exp\left(-C_4\left(\frac{d(x,y)}{t^{1/\beta_1}}\right)^{\frac{\beta_1}{\beta_1-1}}\right),&\text{ if }t<d(x,y),\\
\frac{C_3}{t}\exp\left(-C_4\left(\frac{d(x,y)}{t^{1/\beta_2}}\right)^{\frac{\beta_2}{\beta_2-1}}\right),&\text{ if }t\ge d(x,y).\\
\end{cases}\label{eq_vcsc_grad}
\end{align}
Note that $\partial p_t(x,\cdot)$ exists and has the above estimate $\nu$-a.e. instead of $m$-a.e.. At small scales, $\beta_1=\log15/\log3$, the gradient estimate has the sub-Gaussian behavior corresponding to the Vicsek set, while at large scales, $\beta_2\approx2.09697$, the gradient estimate has the sub-Gaussian behavior corresponding to the Sierpi\'nski carpet. See \cite[Section 3]{BC24} for the results about Lipschitz and gradient estimates for heat kernel on unbounded Vicsek set, that is, the blowup of the Vicsek set by the Vicsek set itself.

\item Let $K$ be the Vicsek set and $M$ the unit interval, see Figure \ref{fig_VC01}.

\begin{figure}[ht]
\begin{center}
\begin{tikzpicture}
\foreach\x in {50pt}{
\foreach\y in {0,1,2,3,4,5}{
\node at (\y*\x,0) {\includegraphics[width=50.5pt]{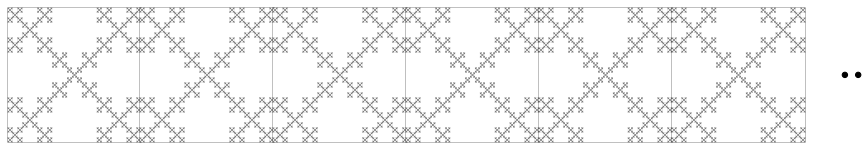}};
}
}
\draw[fill=black] (290pt,0pt) circle (1pt);
\draw[fill=black] (295pt,0pt) circle (1pt);
\draw[fill=black] (300pt,0pt) circle (1pt);
\end{tikzpicture}
\end{center}
\caption{Blow up the Vicsek set by the unit interval}\label{fig_VC01}
\end{figure}
Similar to the blowup of the Vicsek set by the Sierpi\'nski carpet, we also have the Lipschitz estimate, Equation (\ref{eq_vcsc_lip}) and the gradient estimate, Equation (\ref{eq_vcsc_grad}). However, at small scales, $\beta_1=\log15/\log3$, the gradient estimate has the sub-Gaussian behavior corresponding to the Vicsek set, while at large scales, $\beta_2=2$, the gradient estimate has the Gaussian behavior corresponding to the unit interval.
\end{enumerate}

\bibliographystyle{plain}

\end{document}